\documentclass[11pt]{article}
\usepackage{amssymb}
\usepackage{latexsym}
\usepackage{amsfonts}
\usepackage{amsmath}
\newtheorem{thrm}{Theorem}[section]
\newtheorem{lemma}[thrm]{Lemma}
\newtheorem{prop}[thrm]{Proposition}

\newtheorem{cor}[thrm]{Corollary}
\newtheorem{remark}[thrm]{Remark}
\numberwithin{equation}{section}

\usepackage[dvips]{color}

\def\P{\mathbb{P} }
\def\R{\mathbb{R} }
\def\N{\mathbb{N} }
\def\V{\mathbb{V} }
\def\C{\mathcal{C}}
\def\bP{{\bf P} }
\def\D{\mathbb{D}}

\begin{document}
\allowdisplaybreaks
\begin{doublespace}

\title{\Large\bf Central Limit Theorems for Supercritical Superprocesses}
\author{ \bf  Yan-Xia Ren\footnote{The research of this author is supported by NSFC (Grant No.  11271030 and 11128101) and Specialized Research Fund for the
Doctoral Program of Higher Education.\hspace{1mm} } \hspace{1mm}\hspace{1mm}
Renming Song\thanks{Research supported in part by a grant from the Simons
Foundation (208236).} \hspace{1mm}\hspace{1mm} and \hspace{1mm}\hspace{1mm}
Rui Zhang\footnote{Supported by the China Scholarship Council.}
\hspace{1mm} }
\date{}
\maketitle
\begin{abstract}
In this paper, we establish a central limit theorem for a large class of general supercritical
superprocesses with spatially dependent branching mechanisms satisfying a second moment condition.
This central limit theorem generalizes and unifies all the central limit theorems obtained recently in \cite{Mi, RSZ}
for supercritical super Ornstein-Uhlenbeck processes.
The advantage of this central
limit theorem is that it allows us to characterize the limit Gaussian field.
In the case of supercritical super Ornstein-Uhlenbeck processes with non-spatially dependent
branching mechanisms, our central limit theorem reveals more independent structures
of the limit Gaussian field.
\end{abstract}


\medskip


\bigskip

\baselineskip=6.0mm

\section{Introduction}

Central limit theorems for supercritical branching processes were initiated by \cite{KS, KS66}.
In these two papers, Kesten and Stigum established central limit theorems for
supercritical multitype Galton-Watson processes by using the
Jordan canonical form of the expectation matrix $M$.
Then in \cite{Ath69a, Ath69, Ath71}, Athreya proved central limit
theorems for supercritical multi-type continuous time branching processes,
using the Jordan canonical form and the eigenvectors of the matrix $M_t$, the mean matrix at time $t$.
Asmussen and Keiding \cite{AK} used martingale central limit theorems to prove central limit theorems for
supercritical multitype branching processes. In \cite{AH83}, Asmussen and Hering established
spatial central limit theorems for general supercritical branching Markov processes under a certain
condition. However, the condition in \cite{AH83} is not easy to check and essentially the only
examples given in \cite{AH83} of branching Markov processes satisfying this condition are branching
diffusions in bounded smooth domains.
In \cite{RP}, Adamczak and Mi{\l}o\'{s} proved some central limit theorems for supercritical branching
Ornstein-Uhlenbeck processes with binary branching mechanism.
We note that branching Ornstein-Uhlenbeck processes do not satisfy the condition in \cite{AH83}.
In \cite{Mi}, Mi{\l}o\'{s} proved
some central limit theorems for supercritical super Ornstein-Uhlenbeck processes with branching
mechanisms satisfying a fourth moment condition. In \cite{RSZ}, we established central limit
theorems for supercritical super Ornstein-Uhlenbeck processes with (non-spatially dependent) branching
mechanisms satisfying only a second moment condition. More importantly, the central limit
theorems in \cite{RSZ} are more satisfactory since our limit normal random variables are non-degenerate.
In the recent paper \cite{RSZ2}, we obtained central limit theorems for a large class
of general supercritical branching Markov processes with spatially dependent branching mechanisms satisfying
only a second moment condition.
The main results of \cite{RSZ2} are the central limit theorems contained in \cite[Theorems 1.8, 1.9,
1.10, 12]{RSZ2}. \cite[Theorem 1.8]{RSZ2} is the branching Markov process analog of the convergence
of the first and fourth components in Theorem \ref{The:1.3} below.
\cite[Theorem 1.9]{RSZ2} is the branching Markov process analog of Remark \ref{critical} below, while
\cite[Theorems 1.10, 12]{RSZ2} are the branching Markov process analogs of the results in
Remark \ref{large} below.

It is a natural next step to try to establish counterparts of the central limit theorems of \cite{RSZ2} for
general supercritical superprocesses with spatially dependent branching mechanisms satisfying
only a second moment condition.
This is far from trivial.
For a branching Markov process $\{Z_t: t\ge 0\}$, to consider the proper scaling limit of
$\langle f,Z_t\rangle$ as $t\to\infty$,
where $f$ is a test function, it is equivalent to consider the scaling limit of
$\langle f,Z_{{t+s}}\rangle$ as $s\to \infty$ for any $t>0$.
Note that  $Z_{t+s}=\sum_{u\in{\cal L}_t} Z^{u,t}_s$, where ${\cal L}_t$
is the set of particles which are alive at time $t$, and $Z^{u,t}_s$ is the
branching Markov process starting from the particle $u\in \mathcal{L}_t$. So, conditioned on
$Z_t$,  $Z_{t+s}$ is the sum of a finite number of independent terms,
and then basically we only need to consider central limit theorems of independent random variables.
However, a superprocess is an appropriate scaling limit of branching Markov processes,
see \cite{Dawson} and \cite{Li11}, for example. It describes the time evolution of a cloud of
uncountable number of particles, where each particle carries mass $0$ and
moves in space independently according to a Markov process.
The particle picture for superprocesses is not very clear.
Recently \cite{KPR} gave a backbone decomposition of superdiffusions,
where the backbone is a branching diffusion.
One could combine the ideas of \cite{RSZ} with that of \cite{RSZ2} to use the backbone decomposition to prove
central limit theorems for general supercritical superprocesses with
spatial dependent branching mechanisms satisfying only a second moment condition, provided that the
backbone decomposition is known. However, up to now, the backbone decomposition has only been established
for supercritical superdiffusions with spatial dependent branching mechanisms.

In this paper, our assumption on the spatial process is exactly the same as in
\cite{RSZ2}, while our assumptions on the branching mechanism are similar
in spirit to those of \cite{RSZ2}. We will use the excursion measures of the superprocess as a tool to replace the
backbone decomposition. With this new tool, the general methodology of \cite{RSZ2} can be
adapted to the present setting of  general supercritical superprocesses.

Actually, we will go even further in the present paper.
We will prove one central limit
theorem which generalizes and unifies all the central limit theorems of \cite{Mi, RSZ}.
See the Corollaries and Remarks after Theorem \ref{The:1.3}.
The advantage of this central
limit theorem is that it allows us to characterize the limit Gaussian field.
In the case of supercritical super Ornstein-Uhlenbeck processes with non-spatially dependent
branching mechanisms satisfying a second moment condition,
our central limit theorem reveals more independent structures of the limit Gaussian field, see Corollaries
\ref{Cor:1}, \ref{cor:2} and \ref{cor:3}.

\subsection{Spatial process}\label{subs:sp}

Our assumptions on the underlying spatial process are the same as in \cite{RSZ2}.
In this subsection, we recall the assumptions on the spatial process.

$E$ is a locally compact separable metric space and $m$ is a $\sigma$-finite Borel measure
on $E$ with full support.
$\partial$ is a point not contained in $E$ and will be interpreted as the cemetery point.
Every function $f$ on $E$ is automatically extended to $E_{\partial}:=E\cup\{\partial\}$
by setting $f(\partial)=0$.
 We will assume that $\xi=\{\xi_t,\Pi_x\}$ is an $m$-symmetric Hunt process on $E$ and $\zeta:=
 \inf\{t>0: \xi_t=\partial\}$ is the lifetime of $\xi$.
The semigroup of $\xi$ will be denoted by $\{P_t:t\geq 0\}$.
We will always assume
that there exists a family of
continuous strictly positive symmetric functions $\{p_t(x,y):t>0\}$ on $E\times E$ such that
$$
  P_tf(x)=\int_E p_t(x,y)f(y)\,m(dy).
$$
It is well-known that for $p\geq 1$, $\{P_t:t\ge 0\}$ is a
strongly continuous contraction semigroup on $L^p(E, m)$.

Define $\widetilde{a}_t(x):=p_t(x, x)$.
We will always assume that $\widetilde{a}_t(x)$ satisfies the following two conditions:
\begin{description}
  \item[(a)] For any $t>0$, we have
$$
        \int_E \widetilde{a}_t(x)\,m(dx)<\infty.
$$
   \item[(b)] There exists $t_0>0$ such that $\widetilde{a}_{t_0}(x)\in L^2(E,\,m)$.
\end{description}
It is easy to check (see \cite{RSZ2}) that condition $(b)$ above is equivalent to
 \begin{description}
   \item[(b$'$)] There exists $t_0>0$ such that for all $t\ge t_0$, $\widetilde{a}_{t}(x)\in L^2(E,m)$.
 \end{description}

These two conditions are satisfied by a lot of Markov processes. In \cite{RSZ2}, we gave several classes
of examples of Markov processes, including Ornstein-Uhlenbeck processes, satisfying these two conditions.

\subsection{Superprocesses}
In this subsection, we will spell out our assumptions on the superprocess we are going to work with.
Let $\mathcal{B}_b(E)$ ($\mathcal{B}_b^+(E)$) be the set of (positive) bounded Borel measurable functions on $E$.

The superprocess $X=\{X_t:t\ge 0\}$ we are going to work with is determined by three parameters:
a spatial motion $\xi=\{\xi_t, \Pi_x\}$ on $E$ satisfying the
assumptions of the previous subsection,
a branching rate function $\beta(x)$ on $E$ which is a non-negative bounded measurable function and a branching mechanism $\psi$ of the form
\begin{equation}
\psi(x,\lambda)=-a(x)\lambda+b(x)\lambda^2+\int_{(0,+\infty)}(e^{-\lambda y}-1+\lambda y)n(x,dy),
\quad x\in E, \quad\lambda> 0,
\end{equation}
where $a\in \mathcal{B}_b(E)$, $b\in \mathcal{B}_b^+(E)$ and $n$ is a kernel from $E$ to $(0,\infty)$ satisfying
\begin{equation}\label{n:condition}
  \sup_{x\in E}\int_0^\infty y^2 n(x,dy)<\infty.
\end{equation}

Let ${\cal M}_F(E)$
be the space of finite  measures on $E$ equipped with the topology of weak convergence.
The existence of such superprocesses is well-known, see, for instance,
\cite{E.B.} or \cite{Li11}.
$X$ is a  cadlag Markov process taking values in ${\cal M}_F(E)$.
For any $\mu \in \mathcal{M}_F(E)$, we denote the
law of $X$ with initial configuration $\mu$ by $\P_\mu$.
As usual, $\langle f,\mu\rangle:=\int f(x)\mu(dx)$ and $\|\mu\|:=\langle 1,\mu\rangle$.
Then for every
$f\in \mathcal{B}^+_b(E)$ and $\mu \in \mathcal{M}_F(E)$,
\begin{equation}
  -\log \P_\mu\left(e^{-\langle f,X_t\rangle}\right)=\langle u_f(\cdot,t),\mu\rangle,
\end{equation}
where $u_f(x,t)$ is the unique positive solution to the equation
\begin{equation}
  u_f(x,t)+\Pi_x\int_0^t\psi(\xi_s, u_f(\xi_s,t-s))\beta(\xi_s)ds=\Pi_x f(\xi_t),
\end{equation}
where $\psi(\partial,\lambda)=0, \lambda>0$.
Define
\begin{equation}\label{e:alpha}
\alpha(x):=\beta(x)a(x)\quad \mbox{and }
A(x):=\beta(x)\left( 2b(x)+\int_0^\infty y^2 n(x,dy)\right).
\end{equation}
Then, by our assumptions, $\alpha(x)\in \mathcal{B}_b(E)$ and $A(x)\in\mathcal{B}_b(E)$.
Thus there exists $M>0$ such that
\begin{equation}\label{1.5}
  \sup_{x\in E}\left(|\alpha(x)|+A(x)\right)\le M.
\end{equation}
For any $f\in\mathcal{B}_b(E)$ and $(t, x)\in (0, \infty)\times E$, define
\begin{equation}\label{1.26}
   T_tf(x):=\Pi_x \left[e^{\int_0^t\alpha(\xi_s)\,ds}f(\xi_t)\right].
\end{equation}
It is well-known that $T_tf(x)=\P_{\delta_x}\langle f,X_t\rangle$ for every $x\in E$.

It is shown in \cite{RSZ2}
that there exists
a family of continuous strictly positive symmetric functions $\{q_t(x,y),t>0\}$ on $E\times E$ such that
$q_t(x,y)\le e^{Mt}p_t(x,y)$ and
for any $f\in \mathcal{B}_b(E)$,
$$
  T_tf(x)=\int_E q_t(x,y)f(y)\,m(dy).
$$
It follows immediately that, for any $p\ge 1$, $\{T_t: t\ge 0\}$ is a
strongly continuous semigroup on $L^p(E, m)$ and
\begin{equation}\label{Lp}
  \|T_tf\|_p^p\le e^{pMt}\|f\|_p^p.
\end{equation}

Define $a_t(x):=q_t(x, x)$. It follows from the assumptions (a) and (b) in the
previous subsection that $a_t$ enjoys the following properties:
\begin{description}
  \item[(i)] For any $t>0$, we have
$$
        \int_E a_t(x)\,m(dx)<\infty.
$$
  \item[(ii)] There exists $t_0>0$ such that for all $t\ge t_0$, $a_{t}(x)\in L^2(E,m)$.
\end{description}

It follows from (i) above that, for any $t>0$, $T_t$ is  a compact
operator.
The infinitesimal generator $L$ of $\{T_t:t\geq 0\}$ in $L^2(E, m)$ has purely discrete spectrum with eigenvalues
$-\lambda_1>-\lambda_2>-\lambda_3>\cdots$.
It is known that either the number of these eigenvalues is finite, or $\lim_{k\to\infty}\lambda_k=\infty$.
The first eigenvalue $-\lambda_1$ is simple and the eigenfunction
$\phi_1$ associated with $-\lambda_1$ can be chosen to be strictly positive everywhere and continuous.
We will assume that $\|\phi_1\|_2=1$. $\phi_1$ is sometimes denoted as $\phi^{(1)}_1$.
For $k>1$, let $\{\phi^{(k)}_j,j=1,2,\cdots n_k\}$ be
an orthonormal basis of the eigenspace (which is finite dimensional) associated with $-\lambda_k$.
It is well-known that $\{\phi^{(k)}_j,j=1,2,\cdots n_k; k=1,2,\dots\}$ forms a complete orthonormal basis of $L^2(E,m)$
and all the eigenfunctions are continuous.
For any $k\ge 1$, $j=1, \dots, n_k$ and $t>0$, we have $T_t\phi^{(k)}_j(x)=e^{-\lambda_k t}\phi^{(k)}_j(x)$ and
\begin{equation}
\label{1.37}
e^{-\lambda_kt/2}|\phi^{(k)}_j|(x)\le a_t(x)^{1/2}, \qquad x\in E.
\end{equation}
It follows from the relation above that all the eigenfunctions $\phi^{(k)}_j$ belong to $L^4(E, m)$.
For any $x, y\in E$ and $t>0$, we have
\begin{equation}\label{eq:p}
  q_t(x, y)=\sum^\infty_{k=1}e^{-\lambda_k t}\sum^{n_k}_{j=1}\phi^{(k)}_j(x)\phi^{(k)}_j(y),
\end{equation}
where the series is locally uniformly convergent on $E\times E$.
The basic facts recalled in this paragraph are well-known, for instance, one can
refer to  \cite[Section 2]{DS}.

In this paper, we always assume that the superprocess $X$ is supercritical,
that is, $\lambda_1<0$.
Under this assumption, the process $X$ has a strictly positive survival probability,
see the next paragraph.
Note that the number of negative eigenvalues is infinite except in the case
when the total number of eigenvalues is finite.

We will use $\{{\cal F}_t: t\ge0\}$ to denote the filtration of $X$, that is
${\cal F}_t=\sigma(X_s: s\in [0, t])$.
Using the expectation formula of $\langle \phi_1, X_t\rangle$ and the Markov property of $X$,
it is easy to show that (see Lemma \ref{lem:1.2}),
for any nonzero $\mu\in {\cal M}_F(E)$, under $\P_{\mu}$,
the process $W_t:=e^{\lambda_1 t}\langle \phi_1, X_t\rangle$ is a positive martingale.
Therefore it converges:
$$
  W_t \to W_\infty,\quad \P_{\mu}\mbox{-a.s.} \quad \mbox{ as }t\to \infty.
$$
Using the assumption \eqref{n:condition}
we can show that, as $t\to \infty$, $W_t$ also converges in $L^2(\P_{\mu})$,
so $W_\infty$ is non-degenerate and
 the second moment is finite. Moreover, we have $\P_{\mu}(W_\infty)=\langle\phi_1, \mu\rangle$.
Put $\mathcal{E}=\{W_\infty=0\}$, then $\P_{\mu}(\mathcal{E})<1$.
It is clear that  $\mathcal{E}^c\subset\{X_t(E)>0,\forall t\ge 0\}$.

In this paper, we also assume that, for any $t>0$ and $x\in E$,
\begin{equation}\label{extinction}
  \P_{\delta_x}\{\|X_t\|=0\}\in(0,1).
\end{equation}
Here we give a sufficient condition for \eqref{extinction}.
Suppose that $\Phi(z)=\inf_{x\in E}\psi(x,z)\beta(x)$ can be written in the form:
$$
\Phi(z)=\widetilde{a}z+\widetilde{b}z^2+\int^\infty_0(e^{-zy}-1+zy)\widetilde{n}(dy)
$$
with $\widetilde{a}\in \R$, $\widetilde{b}\ge 0$ and $ \widetilde{n}$ is a measure on $(0,\infty)$ satisfying $\int^\infty_0(y\wedge y^2)\widetilde{n}(dy)<\infty$.
If $\widetilde{b}+\widetilde{n}(0,\infty)>0$ and $\Phi(z)$ satisfies
\begin{equation}\label{Phi}
  \int^\infty\frac{1}{\Phi(z)}\,dz<\infty,
\end{equation}
then \eqref{extinction} holds. For the last claim, see, for instance, \cite[Lemma 11.5.1]{Dawson}.

\subsection{Main result}

We will use $\langle\cdot, \cdot\rangle_m$ to denote inner product in
$L^2(E, m)$.
Any $f\in L^2(E,m)$ admits the following eigen-expansion:
\begin{equation}\label{exp:f}
  f(x)=\sum_{k=1}^\infty\sum^{n_k}_{j=1} a_j^k\phi_j^{(k)}(x),
\end{equation}
where $a_j^k=\langle f,\phi_j^{(k)}\rangle_m$ and the series converges in $L^2(E,m)$.
$a^1_1$ will sometimes be written as $a_1$.
For $f\in L^2(E,m)$, define
$$
  \gamma(f):=\inf\{k\geq 1: \mbox{ there exists } j \mbox{ with }
  1\leq j\leq n_k\mbox{ such that }
  a_j^k\neq 0\},
$$
where we use the usual convention $\inf\varnothing=\infty$.

For any $f\in L^2(E,m)$, we define
\begin{eqnarray*}
  f^*(x):=\sum_{j=1}^{n_{\gamma(f)}}a_j^{\gamma(f)}\phi_j^{(\gamma(f))}(x).
\end{eqnarray*}
We note that if $f\in L^2(E,m)$ is nonnegative and $m(x: f(x)>0)>0$,
then
$\langle f,\phi_1\rangle_m>0$
which implies $\gamma(f)=1$
and $f^*(x)=a_1\phi_1(x)=\langle f,\phi_1\rangle_m\phi_1(x).$
The following three subsets of
$L^2(E, m)$
will be needed in the statement of the main result:
$$\C_l:=\left\{g(x)=\sum_{k:\lambda_1>2\lambda_k}\sum^{n_k}_{j=1} b_j^k\phi_j^{(k)}(x): b_j^k\in \R\right\},$$
$$\C_c:=\left\{g(x)=\sum^{n_k}_{j=1} b_j^k\phi_j^{(k)}(x): 2\lambda_k=\lambda_1, b_j^k\in \R\right\}$$
and
$$
\C_s:=\left\{g(x)\in L^2(E,m)\cap L^4(E,m):\lambda_1<2\lambda_{\gamma(g)}\right\}.
$$
Note that $\C_l$ consists of these functions in $L^2(E, m)\cap L^4(E, m)$ that only have nontrivial
projection onto the eigen-spaces corresponding to those ``large'' eigenvalues $-\lambda_k$ satisfying
$\lambda_1>2\lambda_k$. The space $\C_l$ is of finite dimension. The space
$\C_c$ is a subspace  (finite dimensional) of the eigen-space corresponding to
the ``critical'' eigenvalue $-\lambda_k$ with  $\lambda_1=2\lambda_k$. Note that there may not
be a critical eigenvalue and in this case, $\C_c$ is empty. The space
$\C_s$ consists of these functions in $L^2(E, m)\cap L^4(E, m)$ that only have nontrivial
projections onto the eigen-spaces corresponding to those ``small'' eigenvalues $-\lambda_k$ satisfying
$\lambda_1<2\lambda_k$. The space $\C_s$ is of infinite dimensional in general.

In this subsection we give the main result of this paper.
The proof will be given in Section \ref{s:3}.
In the remainder of this paper,
whenever we deal with
an initial configuration
$\mu\in {\cal M}_F(E)$,
we are implicitly assuming that it has compact support.

\subsubsection{Some basic convergence results}

Define
\begin{equation*}
  H_t^{k,j}:=e^{\lambda_k t}\langle\phi_j^{(k)}, X_t\rangle,\quad t\geq 0.
\end{equation*}
Using the same argument as in the proof of \cite[Lemma 3.1]{RSZ2}, we can show that

\begin{lemma}\label{lem:1.2}
$H_t^{k,j}$ is a martingale under $\P_{\mu}$.
Moreover, if $\lambda_1>2\lambda_k$, $\sup_{t>3t_0}\P_{\mu}(H_t^{k,j})^2<\infty$.
Thus the limit
\begin{equation*}
  H_\infty^{k,j}:=\lim_{t\to \infty}H_t^{k.j}
\end{equation*}
exists $\P_{\mu}$-a.s. and in $L^2(\P_{\mu})$.
\end{lemma}

\begin{thrm}\label{The:1.2}
If $f\in L^2(E,m)\cap L^4(E,m)$ with $\lambda_1>2\lambda_{\gamma(f)}$,
then, as $t\to\infty$,
$$
  e^{\lambda_{\gamma(f)}t}\langle f, X_t\rangle\to \sum_{j=1}^{n_{\gamma(f)}}a^{\gamma(f)}_jH^{\gamma(f),j}_\infty,
   \quad \mbox{ in } L^2(\P_{\mu}).
$$
\end{thrm}

{\bf Proof:}
The proof is similar to that of \cite[Theorem 1.6]{RSZ2}. We omit the details here. \hfill$\Box$

\begin{remark}\label{rem:large}
When $\gamma(f)=1$,
$H_t^{1,1}$ reduces to $W_t$, and thus $H_\infty^{1,1}=W_\infty$.
Therefore by Theorem~\ref{The:1.2} and the fact that $a_1=\langle f,\phi_1\rangle_m$,
we get that, as $t\to\infty$,
\begin{equation*}
  e^{\lambda_1 t}\langle f,X_t\rangle\to \langle f,\phi_1\rangle_m W_\infty, \quad \mbox{in } L^2(\P_{\mu}).
\end{equation*}
In particular, the convergence also holds in $\P_{\mu}$-probability.
\end{remark}

\subsubsection {Main Result}

For $f\in \C_s$ and $h\in \C_c$, we define
\begin{equation}\label{e:sigma}
  \sigma_f^2:=\int_0^\infty e^{\lambda_1 s}\langle A(T_s f)^2,\phi_1\rangle_m \,ds
\end{equation}
and
\begin{equation}\label{e:rho}
\rho_h^2:=\left\langle Ah^2,\phi_1\right\rangle_m.
\end{equation}
For $g(x)=\sum_{k: 2\lambda_k<\lambda_1}\sum_{j=1}^{n_k}b_j^k\phi_j^{(k)}(x)\in \C_l,$  we define
\begin{equation}\label{1.61}
I_sg(x):=\sum_{k: 2\lambda_k<\lambda_1}\sum_{j=1}^{n_k}e^{\lambda_k s}b_j^k\phi_j^{(k)} (x)\quad\mbox{and}\quad
 \beta_{g}^2:=\int_0^\infty e^{-\lambda_1 s}\left\langle A(I_sg)^2,\phi_1\right\rangle_m\,ds.
\end{equation}

\begin{thrm}\label{The:1.3}
If $f\in \C_s$, $h\in\C_c$ and
 $g(x)=\sum_{k:2\lambda_k<\lambda_1}\sum_{j=1}^{n_k}b_j^k\phi_j^{(k)}(x)\in\C_l$, then
$ \sigma_f^2<\infty$, $\rho_h^2<\infty$ and $\beta_g^2<\infty$.
Furthermore, it holds that, under $\P_{\mu}(\cdot\mid \mathcal{E}^c)$, as $t\to\infty$,
\begin{eqnarray}\label{result}
   &&\left(e^{\lambda_1 t}\langle \phi_1, X_t\rangle,
 ~\frac{ \langle g,X_t\rangle-\sum_{k:2\lambda_k<\lambda_1}e^{-\lambda_kt}\sum_{j=1}^{n_k}b_j^kH^{k,j}_\infty}
  {\sqrt{\langle \phi_1,X_t\rangle}},
  ~\frac{\langle h , X_t\rangle}{\sqrt{t\langle \phi_1,X_t\rangle}},
  ~\frac{\langle f , X_t\rangle}{\sqrt{\langle \phi_1,X_t\rangle}} \right) \nonumber\\
 && \stackrel{d}{\rightarrow}(W^*,G_3(g),G_2(h),~G_1(f)),
\end{eqnarray}
where $W^*$ has the same distribution as $W_\infty$ conditioned on $\mathcal{E}^c$, $G_3(g)\sim \mathcal{N}(0,\beta_g^2)$, $G_2(h)\sim \mathcal{N}(0,\rho_h^2)$
and $G_1(f)\sim \mathcal{N}(0,\sigma_f^2)$. Moreover, $W^*$, $G_3(g)$, $G_2(h)$ and $G_1(f)$ are independent.
\end{thrm}

This theorem says that, under $\P_{\mu}(\cdot\mid \mathcal{E}^c)$, as $t\to \infty$, the limits
of the second, third and fourth components on the right hand side of \eqref{result} are
nondegenerate normal random variables. Furthermore, the limit normal random variables are
independent. As consequences of this theorem, we could also
get the covariance of the limit random variables
$G_1(f_1)$ and $G_1(f_2)$ when $f_1, f_2\in \C_s$,
the covariance of the limit random variables
$G_2(h_1)$ and $G_2(h_2)$ when $h_1, h_2\in \C_c$,
and the covariance of the limit random variables
$G_3(g_1)$ and $G_3(g_2)$ when $g_1, g_2\in \C_l$.

For $f_1,f_2\in \C_s$, define
$$
\sigma(f_1,f_2)=\int_0^\infty e^{\lambda_1 s}\langle A(T_s f_1)(T_sf_2),\phi_1\rangle_m \,ds.
$$
Note that $\sigma(f,f)=\sigma_f^2$.

\begin{cor}\label{Cor:1}
If $f_1,f_2\in \C_s$,
then, under $\P_{\mu}(\cdot\mid \mathcal{E}^c)$,
$$
 \left(\frac{\langle f_1 , X_t\rangle}{\sqrt{\langle \phi_1,X_t\rangle}},
 \frac{\langle f_2, X_t\rangle}{\sqrt{\langle \phi_1,X_t\rangle}} \right)
 \stackrel{d}{\rightarrow}(G_1(f_1),G_1(f_2)), \quad t\to\infty,
$$
where $(G_1(f_1),G_1(f_2))$ is a bivariate normal random variable with covariance
\begin{equation}\label{sigma(fg)}
 \mbox{\rm Cov}(G_1(f_i),G_1(f_j))=\sigma(f_i,f_j),\quad i,j=1,2.
\end{equation}
Consider the special situation when both the branching mechanism
and the branching rate function are non-spatially dependent,
and $\phi_1$ is a constant function (this is the case of Ornstein-Uhlenbeck processes).
If $f_1=\phi^{(k)}_j$
and $f_2= \phi^{(k')}_{j'}$ are distinct eigenfunctions satisfying $\lambda_1<2\lambda_k$
and  $\lambda_1<2\lambda_{k'}$, then $G_1(f_1)$ and $G_1(f_2)$ are independent.
\end{cor}
\textbf{Proof:}
Using the convergence of the fourth component in Theorem \ref{The:1.3}, we get
\begin{eqnarray*}
  &&{\bP}_\mu \left(\exp\left\{i\theta_1 \frac{\langle f_1,X_t\rangle}{\sqrt{\langle \phi_1,X_t\rangle}}
  +i\theta_2 \frac{\langle f_2,X_t\rangle}{\sqrt{\langle \phi_1,X_t\rangle}}\right\}\mid \mathcal{E}^c\right)\\
  &=&{\bP}_\mu \left(\exp\left\{i\frac{\langle \theta_1 f_1+\theta_2 f_2,X_t\rangle}{\sqrt{\langle \phi_1,X_t\rangle}}\right\} \mid \mathcal{E}^c\right) \\
   &\to& \exp\left\{ -\frac{1}{2}\sigma_{(\theta_1f_1+\theta_2f_2)}^2\right\},\quad \mbox{as}\quad t\to\infty,
\end{eqnarray*}
where
\begin{eqnarray*}
 \sigma_{(\theta_1f_1+\theta_2f_2)}^2 &=& \int_0^\infty e^{\lambda_1 s}\langle A(T_s (\theta_1f_1+\theta_2f_2))^2,\phi_1\rangle_m \,ds \\
   &=&  \theta_1^2\sigma_{f_1}^2+2\theta_1\theta_2\sigma(f_1,f_2)+\theta_2\sigma_{f_2}^2.
\end{eqnarray*}
Note that  $\exp\left\{-\frac{1}{2}\left(\theta_1^2\sigma_{f_1}^2+2\theta_1
\theta_2\sigma(f_1,f_2)+\theta_2\sigma_{f_2}^2\right)\right\} $ is
the characteristic function of   $(G_1(f_1),G_1(f_2))$, which is a  bivariate normal random variable with covariance
$\mbox{\rm Cov}(G_1(f_i),G_1(f_j))=\sigma(f_i,f_j), i,j=1,2$.
The desired result now follows immediately.

In particular, if both the branching mechanism
and the branching rate function are non-spatially dependent, then $A(x)=A$ is a constant.
If $\phi_1$ is a constant function, and $f_1=\phi^{(k)}_j$
and $f_2= \phi^{(k')}_{j'}$ are distinct eigenfunctions satisfying $\lambda_1<2\lambda_k$
and  $\lambda_1<2\lambda_{k'}$, then
$$
\sigma(f_1,f_2)=A\phi_1\int_0^\infty e^{(\lambda_1 -\lambda_k-\lambda_{k'}) s}\langle \phi^{(k)}_j, \phi^{(k')}_{j'}\rangle_m \,ds=0.
$$
and thus  $G_1(f_1)$ and $G_1(f_2)$ are independent.
\hfill$\Box$

For $h_1, h_2\in\C_c$, define
$$
\rho(h_1,h_2)=\langle Ah_1h_2,\phi_1\rangle_m.
$$
Using the convergence of the third component in Theorem \ref{The:1.3}
and an argument similar to that in the proof of Corollary \ref{Cor:1}, we get

\begin{cor}\label{cor:2}
If $h_1, h_2\in\C_c$,
then we have, under $\P_{\mu}(\cdot\mid \mathcal{E}^c)$,
$$
 \left(\frac{\langle h_1 , X_t\rangle}{\sqrt{t\langle \phi_1,X_t\rangle}}, \frac{\langle h_2, X_t\rangle}{\sqrt{t\langle \phi_1,X_t\rangle}} \right)\stackrel{d}{\rightarrow}(G_2(h_1),G_2(h_2)), \quad t\to\infty,
$$
where $(G_2(h_1),G_2(h_2))$ is a bivariate normal random variable with covariance
$$
\mbox{\rm Cov}(G_2(h_i),G_2(h_j))=\rho(h_i,h_j),\quad i,j=1,2.
$$
Consider the special situation when both the branching mechanism
and the branching rate function are non-spatial dependent and $\phi_1$ is a constant function.
If $h_1=\phi^{(k)}_j$
and $h_2= \phi^{(k)}_{j'}$ are distinct eigenfunctions satisfying $\lambda_1=2\lambda_k$,
then $G_2(h_1)$ and $G_2(h_2)$ are independent.
\end{cor}

For $g_1(x)=\sum_{k: 2\lambda_k<\lambda_1}\sum_{j=1}^{n_k}b_j^k\phi_j^{(k)}(x)$ and
$g_2(x)=\sum_{k: 2\lambda_k<\lambda_1}\sum_{j=1}^{n_k}c_j^k\phi_j^{(k)}(x)$,
define
$$
\beta(g_1,g_2)=\int_0^\infty e^{-\lambda_1s}\langle A(I_sg_1)(I_sg_2),\phi_1\rangle_m\,ds.
$$
Using the convergence of the second component in Theorem \ref{The:1.3}
and an argument similar to that in the proof of Corollary \ref{Cor:1}, we get

\begin{cor}\label{cor:3}
If $g_1(x)=\sum_{k: 2\lambda_k<\lambda_1}\sum_{j=1}^{n_k}b_j^k\phi_j^{(k)}(x)$ and
$g_2(x)=\sum_{k: 2\lambda_k<\lambda_1}\sum_{j=1}^{n_k}c_j^k\phi_j^{(k)}(x)$,
then we have, under $\P_{\mu}(\cdot\mid \mathcal{E}^c)$,
\begin{eqnarray*}
&&\left(\frac{\langle g_1,X_t\rangle -\sum_{k: 2\lambda_k<\lambda_1}e^{-\lambda_kt}\sum_{j=1}^{n_k}b_j^kH^{k,j}_\infty}
{\sqrt{\langle \phi_1,X_t\rangle}},
\frac{\langle g_2,X_t\rangle -\sum_{k: 2\lambda_k<\lambda_1}e^{-\lambda_kt}\sum_{j=1}^{n_k}c_j^kH^{k,j}_\infty}
{\sqrt{\langle \phi_1,X_t\rangle}}\right)\\
&&\stackrel{d}{\rightarrow}(G_3(g_1),G_3(g_2)),
\end{eqnarray*}
where $(G_3(g_1),G_3(g_2))$ is a bivariate normal random variable with covariance
$$
\mbox{\rm Cov}(G_3(g_i),G_3(g_j))=\beta(g_i,g_j),\quad i,j=1,2.
$$
Consider the special situation when both the branching mechanism
and the branching rate function are non-spatial dependent and $\phi_1$ is a constant function.
If $g_1=\phi^{(k)}_j$
and $g_2= \phi^{(k')}_{j'}$ are distinct eigenfunctions satisfying $\lambda_1>2\lambda_k$
and  $\lambda_1>2\lambda_{k'}$, then $G_3(g_1)$ and $G_3(g_2)$ are independent.
\end{cor}

\begin{remark}
If $2\lambda_k<\lambda_1$,
then, it holds under $\P_{\mu}(\cdot\mid \mathcal{E}^c)$ that,
as $t\to\infty$,
$$
  \left(e^{\lambda_1 t}\langle \phi_1,X_t\rangle,~\frac{\left(\langle \phi_j^{(k)},X_t\rangle-e^{-\lambda_kt}H^{k,j}_\infty\right)}
  {\langle \phi_1,X_t\rangle^{1/2}} \right)\stackrel{d}{\rightarrow}(W^*,~G_3),
$$
where $G_3\sim\mathcal{N}\left(0,\frac{1}{\lambda_1-2\lambda_k}\langle A(\phi_j^{(k)})^2,\phi_1\rangle_m\right)$.
In particular, for $\phi_1$, we have
\begin{equation*}
  \left(e^{\lambda_1 t}\langle \phi_1,X_t\rangle,
  ~\frac{\left(\langle \phi_1,X_t\rangle-e^{-\lambda_1t}W_\infty\right)}{\langle \phi_1,X_t\rangle^{1/2}}\right)
  \stackrel{d}{\rightarrow}(W^*,~G_3), \quad t\to\infty,
\end{equation*}
where $G_3\sim\mathcal{N}\left(0,-\frac{1}{\lambda_1}\int_E A(x)(\phi_1(x))^3\,m(dx)\right)$.
\end{remark}

All the central limit theorems in \cite{RSZ} are consequences of Theorem \ref{The:1.3}. To see this, we
recall the following notation from \cite{RSZ}.
For $f\in L^2(E,m)$, define
\begin{eqnarray*}
   f_{(s)}(x)& :=& \sum_{k:2\lambda_k<\lambda_1}
  \sum_{j=1}^{n_k}a_j^k\phi_j^{(k)}(x),\\
  f_{(l)}(x)&=& \sum_{k:2\lambda_k>\lambda_1}
  \sum_{j=1}^{n_k}a_j^k\phi_j^{(k)}(x),\\
  f_{(c)}(x)&:=&f(x)-  f_{(s)}(x)-f_{(l)}(x).
\end{eqnarray*}
Then $f_{(s)}\in\C_{l}$, $f_{(c)}\in\C_{c}$ and $f_{(l)}\in\C_{s}$.
Obviously, \cite[Theorem 1.4]{RSZ} is an immediate consequence of the convergence
of the first and fourth components in Theorem \ref{The:1.3}. Now we explain that Theorems 1.6,
1.10 and 1.13 of \cite{RSZ} also follow easily from Theorem \ref{The:1.3}.

\begin{remark}\label{critical}
If $f\in L^2(E,m)\cap L^4(E,m)$ with $\lambda_1=2\lambda_{\gamma(f)}$, then $f=f_{(c)}+f_{(l)}$.
Using the convergence of the fourth component in Theorem \ref{The:1.3} for $f_{(l)}$, it holds under $\P_{\mu}(\cdot\mid \mathcal{E}^c)$ that
$$
\frac{\langle f_{(l)},X_{t}\rangle}{\sqrt{t\langle \phi_1,X_t\rangle}} \stackrel{d}{\to} 0,\quad t\to \infty.
$$
Thus using the convergence of the first and third components in Theorem \ref{The:1.3}, we get, under $\P_{\mu}(\cdot\mid \mathcal{E}^c)$,
$$
  \left(e^{\lambda_1 t}\langle \phi_1, X_t\rangle, ~\frac{\langle f , X_t\rangle}{\sqrt{t\langle \phi_1,X_t\rangle}} \right)\stackrel{d}{\rightarrow}(W^*,~G_2(f_{(c)})), \quad t\to\infty,
$$
where $W^*$ has the same distribution as $W_\infty$ conditioned on $\mathcal{E}^c$
and $G_2(f_{(c)})\sim \mathcal{N}(0,\rho_{f_{(c)}}^2)$. Moreover, $W^*$ and $G_2(f_{(c)})$ are independent.
Thus \cite[Theorem 1.6]{RSZ} is a consequence of Theorem \ref{The:1.3}.
\end{remark}

\begin{remark}\label{large}
Assume $f\in L^2(E,m)\cap L^4(E,m)$ satisfies $\lambda_1>2\lambda_{\gamma(f)}$.

If $f_{(c)}=0$, then $f=f_{(l)}+f_{(s)}$.
Using the convergence of the first, second and fourth components in Theorem \ref{The:1.3},
we get for any nonzero $\mu\in {\cal M}_F(E)$, it holds under $\P_{\mu}(\cdot\mid \mathcal{E}^c)$ that, as $t\to\infty$,
$$
  \left(e^{\lambda_1 t}\langle \phi_1,X_t\rangle,~\frac
  {\left(\langle f,X_t\rangle-\sum_{2\lambda_k<\lambda_1}
  e^{-\lambda_kt}\sum_{j=1}^{n_k}a_j^kH^{k,j}_\infty\right)}{\langle \phi_1,X_t\rangle^{1/2}} \right)
  \stackrel{d}{\rightarrow}(W^*,~G_1(f_{(l)})+G_3(f_{(s)})),
$$
where $W^*$, $G_3(f_{(s)})$ and $G_1(f_{(l)})$ are the same as those in Theorem \ref{The:1.3}.
Since $G_3(f_{(s)})$ and $G_1(f_{(l)})$ are independent,
$G_1(f_{(l)})+G_3(f_{(s)})\sim\mathcal{N}\left(0,\sigma_{f_{(l)}}^2+\beta_{f_{(s)}}^2\right)$.
Thus \cite[Theorem 1.10]{RSZ} is a consequence of Theorem \ref{The:1.3}.

If $f_{(c)}\ne 0$, then as $t\to\infty$,
$$
\frac{\left(\langle f_{(l)}+f_{(s)},X_t\rangle-\sum_{2\lambda_k<\lambda_1}
  e^{-\lambda_kt}\sum_{j=1}^{n_k}a_j^kH^{k,j}_\infty\right)}{\sqrt{t\langle \phi_1,X_t\rangle}}
  \stackrel{d}{\rightarrow} 0.
  $$
Then using the convergence of the first and third components in Theorem \ref{The:1.3}, we get
$$
\left(e^{\lambda_1 t}\langle \phi_1,X_t\rangle,~\frac
  {\left(\langle f,X_t\rangle-\sum_{2\lambda_k<\lambda_1}
  e^{-\lambda_kt}\sum_{j=1}^{n_k}a_j^kH^{k,j}_\infty\right)}{\sqrt{t\langle \phi_1,X_t\rangle}} \right)
  \stackrel{d}{\rightarrow}(W^*,~G_2(f_{(c)})),
$$
where $W^*$ and $G_2(f_{(c)})$ are the same as those in Remark \ref{critical}.
Thus \cite[Theorem 1.13]{RSZ} is a consequence of Theorem \ref{The:1.3}.
\end{remark}

\section{Preliminaries}

\subsection{Excursion measures of $\{X_t,t\ge 0\}$}

We use $\mathbb{D}$ to denote the space of $\mathcal{M}_F({E})$-valued
right continuous functions $t\mapsto \omega_t$ on $(0, \infty)$ having zero as a trap.
We use $(\mathcal{A},\mathcal{A}_t)$ to denote the natural $\sigma$-algebras on
$\mathbb{D}$ generated by the coordinate process.

It is known (see \cite[Section 8.4]{Li11}) that
one can associate with $\{\P_{\delta_x}:x\in E\}$ a family
of $\sigma$-finite measures $\{\N_x:x\in E\}$ defined on $(\D,\mathcal{A})$ such  that $\mathbb{N}_x(\{0\})=0$,
\begin{equation}
 \int_{\D}(1- e^{-\langle f, \omega_{t} \rangle})\mathbb{N}_x(d\omega)
= -\log \mathbb{P}_{\delta_x}(e^{-\langle f, X_{t} \rangle}) ,
\quad f\in {\cal B}^+_b(E),\ t> 0,
\label{DK}
\end{equation}
and, for every $0<t_1<\cdots<t_n<\infty$, and nonzero $\mu_1,\cdots,\mu_n\in M_F(E)$,
\begin{eqnarray}\label{TN}
  &&\mathbb{N}_x(\omega_{t_1}\in d\mu_1,\cdots,\omega_{t_n}\in d\mu_n) \nonumber\\
 && = \mathbb{N}_x(\omega_{t_1}\in d\mu_1)\P_{\mu_1}(X_{t_2-t_1}\in d\mu_2)\cdots \P_{\mu_{n-1}}(X_{t_n-t_{n-1}}\in d\mu_n).
\end{eqnarray}
For earlier work on excursion measures of superprocesses, see \cite{elk-roe, Li03, E.B2.}.

For any $\mu\in M_F(E)$, let $N(d\omega)$
be a Poisson random measure on the space $\mathbb{D}$
with intensity $\int_ E\mathbb{N}_x(d\omega)\mu(dx)$, in a probability space $(\widetilde\Omega, \widetilde{\cal F}, \mathbf{P}_\mu)$.
Define another process $\{\Lambda_t: t\ge 0\}$ by $\Lambda_0=\mu$ and
$$\Lambda_t:=\int_{\mathbb{D}}\omega_tN(d\omega),\quad t>0.$$
Let $\widetilde{\mathcal{F}}_t$ be the $\sigma$-algebra generated by the random variables $\{N(A):A\in \mathcal{A}_t\}$.
Then, $\{\Lambda,(\widetilde{\mathcal{F}}_t)_{t\ge 0}, \mathbf{P}_{\mu}\}$ has the same law as $\{X,(\mathcal{F}_t)_{t\ge 0},\P_{\mu}\}$, see \cite[Theorem 8.24]{Li11} for a proof.

Now we list some properties of $\mathbb{N}_x$. The proofs are similar to those in \cite[Corollary 1.2, Proposition 1.1]{E.B2.}.
\begin{prop}\label{Moments of N}
If $\mathbb{P}_{\delta_x}|\langle f, X_t\rangle|<\infty,$ then
\begin{equation}\label{N1}
   \int_{\D}\langle f, \omega_t\rangle\,\N_x(d\omega)=\mathbb{P}_{\delta_x}\langle f, X_t\rangle.
\end{equation}
If $\mathbb{P}_{\delta_x}\langle f, X_t\rangle^2<\infty,$ then
\begin{equation}\label{N2}
  \int_{\D}\langle f, \omega_t\rangle^2\,\N_x(d\omega)=\mathbb{V}ar_{\delta_x}\langle f, X_t\rangle.
\end{equation}
\end{prop}
\begin{prop}\label{neq0}
\begin{equation}\label{N3}
  \mathbb{N}_x(\|\omega_t\|\neq 0)=-\log\mathbb{P}_{\delta_x}(\|X_t\|=0).
\end{equation}
\end{prop}
\begin{remark}\label{Rek5}
By \eqref{extinction} and Proposition \ref{neq0}, for each $t>0$ and $x\in E$, we have
$$
0<\mathbb{N}_x(\|\omega_t\|\neq 0)<\infty.
$$
Thus, we can define another probability measure $\widetilde{\N}_x$ on $\mathbb{D}$ as follows:
\begin{equation}\label{tilde N}
  \widetilde{\N}_x(B)=\frac{\mathbb{N}_x
  \left(B\cap\{\|\omega_{1}\|\neq 0\}\right)}{\mathbb{N}_x(\|\omega_{1}\|\neq 0)}.
\end{equation}
\end{remark}

Notice that, for $f\in L^2(E,m)$, $\N_x(\langle |f|,\omega_t\rangle)=T_t|f|(x)<\infty$, which implies that $\N_x(\langle |f|,\omega_t\rangle=\infty)=0$.
Thus, for $f\in L^2(E,m)$,
\begin{eqnarray*}
  \P_{\mu}\left(e^{i\theta\langle f,X_t\rangle}\right)&=&\mathbf{P}_{\mu}\left(e^{i\theta\langle f,\Lambda_t\rangle}\right)
  =\mathbf{P}_{\mu}\left(e^{i\theta\int_{\D}\langle f,\omega_t\rangle \,N(d\omega)}\right)\\
   &=& \exp\left\{\int_E\int_\D \left(e^{i\theta\langle f,\omega_t\rangle}-1\right)\N_x(d\omega)\mu(dx)\right\}.
\end{eqnarray*}
Thus, by the Markov property of superprocesses, we have
\begin{equation}\label{cf}
  \mathbb{P}_\mu\left[\exp\left\{i\theta\langle f, X_{t+s}\rangle\right\}|X_t\right]=\P_{X_t}\left(e^{i\theta\langle f,X_s\rangle}\right)
  =\exp\left\{\int_{E}\int_{\mathbb{D}}(e^{i\theta\langle f, \omega_s\rangle}-1)\mathbb{N}_x(d\omega)X_t(dx)\right\}.
\end{equation}

\subsection{Estimates on the moments of  $X$}
In the remainder of this paper we will use the  following notation:
for two positive functions $f$ and $g$ on $E$, $f(x)\lesssim g(x)$ means that there exists a constant $c>0$ such that
$f(x)\le cg(x)$ for all $x\in E$.

First, we recall some results about the semigroup $(T_t)$, the proofs of which can be found in \cite{RSZ2}.

\begin{lemma}\label{lem:expansion}
 For any $f\in L^2(E,m)$, $x\in E$ and $t>0$, we have
\begin{equation}\label{1.17}
T_t f(x)=\sum_{k=\gamma(f)}^\infty e^{-\lambda_k t}\sum_{j=1}^{n_k}a_j^k\phi^{(k)}_j(x)
\end{equation}
and
\begin{equation}\label{1.25}
  \lim_{t\to\infty}e^{\lambda_{\gamma(f)}t}T_t f(x)=\sum_{j=1}^{n_{\gamma(f)}}a_j^{\gamma(f)}\phi_j^{(\gamma(f))}(x),
\end{equation}
where the series in \eqref{1.17} converges absolutely and uniformly in any compact subset of $E$.
Moreover, for any $t_1>0$,
\begin{eqnarray}
   &&\sup_{t>t_1}e^{\lambda_{\gamma(f)}t}|T_t f(x)|\le e^{\lambda_{\gamma(f)}t_1}\|f\|_2\left(\int_E a_{t_1/2}(x)\,m(dx)\right)a_{t_1}(x)^{1/2},\label{1.36}\\
  &&\sup_{t>t_1}e^{(\lambda_{\gamma(f)+1}-\lambda_{\gamma(f)})t}
  \left|e^{\lambda_{\gamma(f)}t}T_t f(x)-f^*(x)\right|
  \le e^{\lambda_{\gamma(f)+1}t_1}\|f\|_2\left(\int_E a_{t_1/2}(x)\,m(dx)\right)(a_{t_1}(x))^{1/2}.\nonumber\\
  &&\label{1.43}
\end{eqnarray}
\end{lemma}

\begin{lemma}\label{lem:rsnew}
Suppose that $\{f_t(x):t>0\}$ is a family of functions in $L^2(E, m)$.
 If $\lim_{t\to\infty} \|f_t\|_2=0$, then for any $x\in E$,
$$
\lim_{t\to\infty}e^{\lambda_1t}T_t f_t(x)=0.
$$
\end{lemma}

Recall the second moments of the superprocess $\{X_t: t\ge 0\}$ (see, for example, \cite[Corollary 2.39]{Li11}):
for $f\in \mathcal{B}_b(E)$, we have for any $t >0$,
\begin{equation}\label{1.9}
  \P_{\mu}\langle f,X_t\rangle^2=\left(\P_{\mu}\langle f,X_t\rangle\right)^2+\int_E\int_0^tT_{s}[A(T_{t-s}f)^2](x)\,ds\mu(dx).
\end{equation}
Thus,
\begin{equation}\label{1.13}
   {\V}{\rm ar}_{\mu}\langle f,X_t\rangle=\langle{\V}{\rm ar}_{\delta_\cdot}\langle f,X_t\rangle, \mu\rangle=\int_E\int_0^tT_{s}[A(T_{t-s}f)^2](x)\,ds\mu(dx),
\end{equation}
where $\mathbb{V}{\rm ar}_{\mu}$ stands for the variance under $\P_{\mu}$.
Note that the second moment formula \eqref{1.9} for superprocesses
is different from  that of \cite[(2.11)]{RSZ2} for branching Markov processes.

For any $f\in L^2(E,m)\cap L^4(E,m)$ and $x\in E$, since $(T_{t-s}f)^2(x)\le e^{M(t-s)}T_{t-s}(f^2)(x)$, we have
\begin{equation*}
  \int_0^tT_{s}[A(T_{t-s}f)^2](x)\,ds\le e^{Mt}T_t(f^2)(x)<\infty.
\end{equation*}
Thus, using a routine limit argument, one can easily check
that \eqref{1.9} and \eqref{1.13} also hold for $f\in L^2(E,m)\cap L^4(E,m)$.

\begin{lemma} \label{lem:2.2}
Assume that $f\in L^2(E,m)\cap L^4(E,m)$.
\begin{description}
  \item{(1)}
   If $\lambda_1<2\lambda_{\gamma(f)}$, then for any $x\in E$,
    \begin{equation}\label{limit-mean}
   \lim_{t\to \infty}e^{\lambda_1t/2}\P_{\delta_x}\langle f, X_t\rangle=0,
   \end{equation}
  \begin{equation}\label{2.2}
   \lim_{t\to \infty}e^{\lambda_1t}{\V}{\rm ar}_{\delta_x}\langle f,X_t\rangle
  =\sigma^2_f \phi_1(x),
  \end{equation}
  where $\sigma^2(f)$ is defined by \eqref{e:sigma}.
  Moreover, for $(t, x)\in (3t_0, \infty)\times E$, we have
  \begin{equation}\label{2.10}
     e^{\lambda_1t}{\V}{\rm ar}_{\delta_x}\langle f,X_t\rangle\lesssim a_{t_0}(x)^{1/2}.
  \end{equation}

  \item{(2)}
  If $\lambda_1=2\lambda_{\gamma(f)}$, then for any $(t, x)\in (3t_0, \infty)\times E$,
  \begin{equation}\label{1.49}
   \left|t^{-1}e^{\lambda_1 t}{\V}{\rm ar}_{\delta_x}\langle f,X_t\rangle-\rho_{f^*}^2\phi_1(x)\right|
   \lesssim t^{-1}a_{t_0}(x)^{1/2},
  \end{equation}
  where $\rho^2_{f^*}$ is defined by \eqref{e:rho}.

  \item{(3)}
If $\lambda_1>2\lambda_{\gamma(f)}$, then for any $x\in E$,
  \begin{equation}\label{1.23}
    \lim_{t\to \infty}e^{2\lambda_{\gamma(f)}t}{\V}{\rm ar}_{\delta_x}\langle f,X_t\rangle=\eta_f^2(x),
  \end{equation}
  where
  \begin{equation*}\label{e:etaf}
   \eta_f^2(x):=\int_0^\infty e^{2\lambda_{\gamma(f)}s}T_s(A(f^*)^2)(x)\,ds.
  \end{equation*}
  Moreover,
  for any $(t, x)\in (3t_0, \infty)\times E$,
  \begin{equation}\label{1.60}
   e^{2\lambda_{\gamma(f)}t}\P_{\delta_x}\langle f,X_t\rangle^2\lesssim a_{t_0}(x)^{1/2}.
  \end{equation}
\end{description}
\end{lemma}
\textbf{Proof:}
Since the first moment formulas for superprocesses and branching Markov processes
are the same, we get \eqref{limit-mean} easily. Although the second moment formula
for superprocesses is different from that for branching Markov processes,
we can still get all results on the variance of the superprocess $X$ from the proof of \cite[Lemma 2.3]{RSZ2}. In fact,
$$
 {\V}{\rm ar}_{x}\langle f,X_t\rangle=\int_0^tT_{s}[A(T_{t-s}f)^2](x)\,ds.
$$
The limit behaviour of the right side of the above equation, as $t\to\infty$, was given in the proof of \cite[Lemma 2.3]{RSZ2}
\hfill$\Box$

\begin{lemma}\label{lem:nnew}
Assume that $f\in L^2(E,m)\cap L^4(E,m)$. If $\lambda_1<2\lambda_{\gamma(f)}$, then for any $(t, x)\in (3t_0, \infty)\times E$,
\begin{equation}\label{small:new}
   \left|e^{\lambda_1t}{\V}{\rm ar}_{\delta_x}\langle f,X_t\rangle
   - \sigma^2_f\phi_1(x)\right|\lesssim \left(e^{(\lambda_1-2\lambda_{\gamma(f)})t}+e^{(\lambda_1-\lambda_2)t}\right)a_{t_0}(x)^{1/2}.
  \end{equation}
\end{lemma}
\textbf{Proof:}
By \eqref{1.13}, we get, for $t>3t_0$,
\begin{eqnarray}
  &&\left|e^{\lambda_1t}{\V}{\rm ar}_{\delta_x}\langle f,X_t\rangle
 - \int_0^\infty e^{\lambda_1 s}\langle A(T_s f)^2,\phi_1\rangle_m \,ds \phi_1(x)\right| \nonumber\\
  &=& \left|e^{\lambda_1t}\int^t_0T_{t-s}[A(T_sf)^2](x)\,ds-\int_0^\infty e^{\lambda_1 s}
 \langle A(T_s f)^2,\phi_1\rangle_m \,ds \phi_1(x)\right|\nonumber \\
  &\le& e^{\lambda_1t}\int^{t-t_0}_0\left|T_{t-s}[A(T_sf)^2](x)-e^{-\lambda_1(t-s)}
 \langle A(T_s f)^2,\phi_1\rangle_m\phi_1(x)\right|\,ds\nonumber\\
  &&+e^{\lambda_1t}\int_{t-t_0}^tT_{t-s}[A(T_sf)^2](x)\,ds+\int_{t-t_0}^\infty e^{\lambda_1 s}
  \langle A(T_s f)^2,\phi_1\rangle_m \,ds \phi_1(x)\nonumber\\
  &=:& V_1(t,x)+V_2(t,x)+V_3(t,x).
\end{eqnarray}
For $V_2(t,x)$, by \cite[(2.26)]{RSZ2}, we have
\begin{equation}\label{1.14}
  V_2(t,x)\lesssim e^{(\lambda_1-2\lambda_{\gamma(f)})t}a_{t_0}(x)^{1/2}.
\end{equation}
For $V_3(t,x)$, by \eqref{1.36}, for $s>t-t_0>t_0$, $|T_sf(x)|\lesssim e^{-\lambda_{\gamma(f)}s}a_{t_0}(x)^{1/2}$.
By \eqref{1.37}, $\phi_1(x)\le e^{\lambda_1t_0/2}a_{t_0}(x)^{1/2}$. Thus, we get
\begin{eqnarray}\label{1.15}
  V_3(t,x) &\lesssim& \int_{t-t_0}^\infty e^{(\lambda_1-2\lambda_{\gamma(f)}) s} \,ds
 \langle a_{t_0},\phi_1\rangle_m \phi_1(x) \nonumber\\
   &\lesssim &  e^{(\lambda_1-2\lambda_{\gamma(f)})t}a_{t_0}(x)^{1/2}.
\end{eqnarray}
Finally, we consider $V_1(t,x)$.
Using \eqref{1.43} with $f$ replaced by $g:=A(T_sf)^2$ and noticing that $\gamma(g)=1$ and $g^*(x)=\langle A(T_s f)^2,\phi_1\rangle_m\phi_1(x)$, for $t-s>t_0$, we have
\begin{eqnarray*}
  \left|T_{t-s}[A(T_sf)^2](x)-e^{-\lambda_1(t-s)}
  \langle A(T_s f)^2,\phi_1\rangle_m\phi_1(x)\right|
  &\lesssim& e^{-\lambda_2(t-s)}\|A(T_s f)^2\|_2a_{t_0}(x)^{1/2}.
\end{eqnarray*}
For $s>t_0$, by \eqref{1.36}, $|T_sf(x)|\lesssim e^{-\lambda_{\gamma(f)}s}a_{t_0}(x)^{1/2}$.
Thus,
$$
\|A(T_s f)^2\|_2\lesssim e^{-2\lambda_{\gamma(f)}s}\|a_{t_0}\|_2.
$$
For $s\le t_0$, by \eqref{Lp}, it is easy to get
$$
\|A(T_s f)^2\|_2\le M\|T_s f\|_4^2\le Me^{2Ms}\|f\|_4^2.
$$
Therefore, we have
\begin{eqnarray}\label{1.16}
  V_1(t,x) &\lesssim&
  e^{\lambda_1t}\int_{t_0}^{t-t_0} e^{-\lambda_2(t-s)}e^{-2\lambda_{\gamma(f)}s}\,ds\,a_{t_0}(x)^{1/2}
  +e^{\lambda_1t}\int_{0}^{t_0}e^{-\lambda_2(t-s)}\,ds\,a_{t_0}(x)^{1/2}\nonumber\\
   &\lesssim& \left(e^{(\lambda_1-2\lambda_{\gamma(f)})t}+e^{(\lambda_1-\lambda_2)t}\right)a_{t_0}(x)^{1/2}.
\end{eqnarray}
Now \eqref{small:new} follows immediately from \eqref{1.14}, \eqref{1.15} and \eqref{1.16}.\hfill$\Box$

\begin{lemma}\label{small and critical}
Assume that $f\in L^2(E,m)\cap L^4(E,m)$ with $\lambda_1<2\lambda_{\gamma(f)}$
and $h\in L^2(E,m)\cap L^4(E,m)$ with $\lambda_1=2\lambda_{\gamma(h)}$. Then, for any $(t, x)\in (3t_0, \infty)\times E$,
\begin{equation}\label{cov:sc}
 \mathbb{C}{\rm ov}_{\delta_x}(e^{\lambda_1t/2}\langle f,X_t\rangle, t^{-1/2}e^{\lambda_1t/2}\langle h,X_t\rangle)\lesssim t^{-1/2}(a_{t_0}(x))^{1/2},
\end{equation}
where $\mathbb{C}{\rm ov}_{\delta_x}$ is the covariance under $\P_{\delta_x}$.
\end{lemma}
\textbf{Proof:}
By \eqref{1.13}, we have
\begin{eqnarray*}
  &&\left|\mathbb{C}{\rm ov}_{\delta_x}(e^{\lambda_1t/2}\langle f,X_t\rangle, t^{-1/2}e^{\lambda_1t/2}\langle h,X_t\rangle)\right| \\
  &=& t^{-1/2}e^{\lambda_1t}\frac{1}{4}\left|\left(\mathbb{V}{\rm ar}_{\delta_x}\langle (f+h),X_t\rangle- \mathbb{V}{\rm ar}_{\delta_x}\langle (f-h),X_t\rangle\right)\right|\\
   &=& t^{-1/2}e^{\lambda_1t}\left|\int_0^t T_{t-s}\left[A(T_sf)(T_sh)\right](x)\,ds\right|\\
   &\le& t^{-1/2}e^{\lambda_1t}\left(\int_0^{t-t_0} T_{t-s}[A\left|(T_sf)(T_sh)\right|](x)\,ds
   +\int_{t-t_0}^t T_{t-s}[A\left|(T_sf)(T_sh)\right|](x)\,ds\right)\\
   &=:&V_4(t,x)+V_5(t,x).
\end{eqnarray*}
First, we deal with $V_4(t,x)$. By \eqref{1.36}, for $t-s>t_0$,
$$T_{t-s}[A\left|(T_sf)(T_sh)\right|](x)\lesssim e^{-\lambda_1(t-s)}\|A(T_sf)(T_sh)\|_2(a_{t_0}(x))^{1/2}.$$
If $s>t_0$, then by \eqref{1.36}, we get
$$
\|A(T_sf)(T_sh)\|_2\lesssim e^{-(\lambda_1/2+\lambda_{\gamma(f)})s}\|a_{t_0}\|_2.
$$
If $s\le t_0$, by \eqref{Lp}, it is easy to get
$$
\|A(T_sf)(T_sh)\|_2\le M\|T_sf\|_4\|T_sh\|_4\le Me^{2Ms}\|f\|_4\|h\|_4.
$$
Therefore, we have
\begin{eqnarray}\label{1.18}
  V_4(t,x)&\lesssim&
  t^{-1/2}e^{\lambda_1t}\left(\int_{t_0}^{t-t_0} e^{-\lambda_1(t-s)}e^{-(\lambda_1/2+\lambda_{\gamma(f)})s}\,ds
  +\int_{0}^{t_0}e^{-\lambda_1(t-s)}\,ds\right)\,a_{t_0}(x)^{1/2}\nonumber\\
  &=& t^{-1/2}\left(\int_{t_0}^{t-t_0}e^{(\lambda_1/2-\lambda_{\gamma(f)})s}\,ds
  +\int_{0}^{t_0}e^{\lambda_1s}\,ds\right)\,a_{t_0}(x)^{1/2}\nonumber\\
   &\lesssim& t^{-1/2}a_{t_0}(x)^{1/2}.
\end{eqnarray}

For $V_5(t,x)$, if $s>t-t_0\ge 2t_0$, then by\eqref{1.36}, we get
\begin{eqnarray}\label{1.20}
  V_5(t,x)&\lesssim& t^{-1/2}e^{\lambda_1t}\int_{t-t_0}^t e^{-(\lambda_1/2+\lambda_{\gamma(f)})s}T_{t-s}(a_{2t_0})(x)\,ds\nonumber\\
   &=&  t^{-1/2}e^{(\lambda_1/2-\lambda_{\gamma(f)})t}\int_0^{t_0}e^{(\lambda_1/2+\lambda_{\gamma(f)})s}T_{s}(a_{2t_0})(x)\,ds\nonumber\\
   &\lesssim& t^{-1/2}e^{(\lambda_1/2-\lambda_{\gamma(f)})t}\int_0^{t_0}T_{s}(a_{2t_0})(x)\,ds\nonumber\\
   &\lesssim&t^{-1/2}(a_{t_0}(x))^{1/2}.
\end{eqnarray}
The last inequality follows from the fact that
\begin{equation}\label{1.19}
  \int_{0}^{t_0}T_s(a_{2t_0})(x)\,ds\lesssim a_{t_0}(x)^{1/2},
\end{equation}
which is \cite[(2.25)]{RSZ2}.
Therefore, by \eqref{1.18} and \eqref{1.20}, we get \eqref{cov:sc} immediately.\hfill$\Box$

\section{Proof of the main theorem}\label{s:3}

In this section, we will prove the main result of this paper.
The general methodology is similar to that of \cite{RSZ2}, the difference being
that we use the excursion measures of the superprocess rather than the backbone
decomposition (which is not yet available in the general setup of this paper)
of superprocess.

We first recall some facts about weak convergence which will be used later.
For $f:\mathbb{R}^n\to\mathbb{R}$, let $\|f\|_L:=\sup_{x\ne y}|f(x)-f(y)|/\|x-y\|$ and
  $\|f\|_{BL}:=\|f\|_{\infty}+\|f\|_L$. For any distributions $\nu_1$ and $\nu_2$ on $\mathbb{R}^n$, define
\begin{equation*}
  d(\nu_1,\nu_2):=\sup\left\{\left|\int f\,d\nu_1-\int f\,d\nu_2\right|~:~\|f\|_{BL}\leq1\right\}.
\end{equation*}
Then $d$ is a metric.
It follows from \cite[Theorem 11.3.3]{Dudley} that the topology generated by $d$
is equivalent to the weak convergence topology.
 From the definition, we can easily see that, if $\nu_1$ and $\nu_2$ are the distributions of two $\R^n$-valued random variables $X$ and $Y$ respectively, then
\begin{equation}\label{5.20}
  d(\nu_1,\nu_2)\leq E\|X-Y\|\leq\sqrt{ E\|X-Y\|^2}.
\end{equation}

The following simple fact will be used several times later in this section:
\begin{equation}\label{3.20}
  \left|e^{ix}-\sum_{m=0}^n\frac{(ix)^m}{m!}\right|\leq \min\left(\frac{|x|^{n+1}}{(n+1)!}, \frac{2|x|^n}{n!}\right).
\end{equation}

Before we prove Theorem \ref{The:1.3}, we prove several lemmas first.
The first lemma below says that the result in Lemma \ref{lem:1.2} also holds under ${\N}_x$.
Recall the probability measure $\widetilde{\N}_x$ defined in \eqref{tilde N}.
On the measurable space $\left(\D, \mathcal{A}\right)$, define
$$
  \widetilde{H}_t^{k,j}(\omega):=e^{\lambda_kt}\langle\phi_j^{(k)}, \omega_t\rangle,\quad t\geq 0,\quad \omega\in \D.
$$
\begin{lemma}\label{lem:2.5}
For $x\in E$,
if $\lambda_1>2\lambda_k$, then
the limit
\begin{equation*}
 \widetilde{ H}_\infty^{k,j}:=\lim_{t\to \infty}\widetilde{H}_t^{k,j}
\end{equation*}
exists $\mathbb{N}_x$-a.e., in $L^1(\mathbb{N}_x)$ and in $L^2(\mathbb{N}_x)$.
\end{lemma}

{\bf Proof:}
On the set $\{\omega\in \D:\|\omega_{1}\|=0 \}$, we have $\omega_t=0$, $t>1$, thus, $\widetilde{H}_\infty^{k,j}(\omega)=0$. Thus, we only need to show $\widetilde{H}_\infty^{k,j}$ exists $\widetilde{\N}_x$-a.s. and in $L^2(\widetilde{\N}_x)$.

For $t>s\ge 1$, since $\{\|\omega_{1}\|= 0\}\subset \{\|\omega_{s}\|= 0 \}\subset \{\|\omega_{t}\|=0 \}$, we have
$$\mathbb{N}_x\left(\langle \phi_j^{(k)},\omega_t\rangle; \|\omega_{1}\|\neq 0|\mathcal{A}_s\right)=\mathbb{N}_x\left(\langle \phi_j^{(k)},\omega_t\rangle|\mathcal{A}_s\right)=\P_{\omega_s}\left(\langle \phi_j^{(k)},X_{t-s}\rangle\right)=e^{-\lambda_k(t-s)}\langle \phi_j^{(k)},\omega_s\rangle,$$
which implies $\{\widetilde{H}_t^{k,j},t\ge 1\}$ is a martingale under $\widetilde{\N}_x$.
By \eqref{N2}, we have
$$\mathbb{N}_x\left(\langle \phi_j^{(k)},\omega_t\rangle^2; \|\omega_{1}\|\neq 0\right)=\mathbb{N}_x\left(\langle \phi_j^{(k)},\omega_t\rangle^2\right)=\mathbb{V}ar_{\delta_x}\langle \phi_j^{(k)},X_t\rangle.$$
Then by Lemma \ref{lem:1.2}, we easily get $\limsup_{t\to\infty}\widetilde{\N}_x(\widetilde{H}_t^{k,j})^2<\infty$,
which implies $\widetilde{H}_\infty^{k,j}$ exists $\widetilde{\N}_x$-a.s. and in $L^2(\widetilde{\N}_x)$. \hfill$\Box$

\begin{lemma}\label{lem:small}
If $f\in \C_s$,
then $ \sigma_f^2<\infty$ and,
for any nonzero $\mu\in {\cal M}_F(E)$, it holds under $\P_{\mu}$ that
$$
  \left(e^{\lambda_1 t}\langle \phi_1, X_t\rangle, ~e^{\lambda_1t/2}\langle f , X_t\rangle \right)\stackrel{d}{\rightarrow}\left(W_\infty,~G_1(f)\sqrt{W_\infty}\right), \quad t\to\infty,
$$
where $G_1(f)\sim \mathcal{N}(0,\sigma_f^2)$. Moreover, $W_\infty$ and $G_1(f)$ are independent.
\end{lemma}

\textbf{Proof:}\quad
We need to consider the limit of the ${\mathbb R}^2$-valued random variable $U_1(t)$ defined by
\begin{equation}\label{5.16}
   U_1(t):=\left(e^{\lambda_1 t}\langle \phi_1,X_t\rangle,\, e^{\lambda_1t/2}\langle f, X_t\rangle\right),
\end{equation}
or equivalently, we need to consider the limit of $U_1(t+s)$ as $t\to \infty$  for any $s>0$.
The main idea is as
follows. For $s,t> t_0$,
\begin{equation}\label{decom-U1}\begin{array}{rll}U_1(s+t)&=&\displaystyle\left(e^{\lambda_1 (t+s)}\langle \phi_1,X_{t+s}\rangle,\, e^{\lambda_1(t+s)/2}\langle f,X_{t+s}\rangle-e^{\lambda_1(t+s)/2}\langle T_s f,X_{t}\rangle\right)\\
&&\displaystyle+\left(0,\, e^{\lambda_1(t+s)/2}\langle T_s f,X_{t}\rangle\right).\end{array}\end{equation}
The double limit, first as $t\to\infty$ and then $s\to \infty$, of the first term of the right side of \eqref{decom-U1}
is equal to the double limit, first as $t\to\infty$ and then $s\to \infty$, of another
${\mathbb R}^2$-valued random variable $U_2(s,t)$
where \begin{eqnarray*}
U_2(s,t):=\left(e^{\lambda_1 t}\langle \phi_1,X_t\rangle,e^{\lambda_1(t+s)/2}\langle f,X_{t+s}\rangle-e^{\lambda_1(t+s)/2}\langle T_s f,X_{t}\rangle\right).
\end{eqnarray*}
We will prove that the second term on the right hand side of \eqref{decom-U1} has no contribution to the
double limit, first as $t\to\infty$ and then $s\to \infty$, of the left hand side
(see, \eqref{limsuplimsup} below).

We claim that, under $\mathbb{P}_\mu$,
\begin{equation}\label{10.5}
  U_2(s,t)\stackrel{d}{\to}\left(W_\infty, \sqrt{W_\infty}G_1(s)\right), \quad \mbox{ as } t\to\infty,
\end{equation}
where $G_1(s)\sim\mathcal{N}(0,\sigma^2_f(s))$ with $\sigma^2_f(s)$ to be given later.
In fact, denote the characteristic function of $U_2(s,t)$ under $\mathbb{P}_\mu$ by
$\kappa(\theta_1,\theta_2,s,t)$:
\begin{eqnarray}
 && \kappa(\theta_1,\theta_2,s,t)\nonumber\\
 &=&\mathbb{P}_{\mu}\left(\exp\left\{i\theta_1e^{\lambda_1 t}\langle \phi_1,X_t\rangle
 +i\theta_2e^{\lambda_1(t+s)/2}\langle f,X_{t+s}\rangle-i\theta_2e^{\lambda_1(t+s)/2}\langle T_s f,X_{t}\rangle\right\}\right)\nonumber\\
  &=&\mathbb{P}_{\mu}\left(\exp\left\{i\theta_1e^{\lambda_1 t}\langle \phi_1,X_t\rangle\right.\right.\nonumber\\
  &&\left.\left.+\int_{E}\int_{\mathbb{D}}
     (\exp\left\{ i\theta_2e^{\lambda_1(t+s)/2}\langle f, \omega_s\rangle\right\}
  -1-i\theta_2e^{\lambda_1(t+s)/2}\langle f,\omega_s\rangle)\mathbb{N}_x(d\omega)X_t(dx)\right\}\right)
  ,\label{10.9}
\end{eqnarray}
where in the last equality we used the Markov property of $X$, \eqref{N1} and \eqref{cf}.
Define
$$R_s(\theta,x)=\int_{\mathbb{D}}\left(\exp\{\langle i\theta f,\omega_s\rangle\}-1-i\theta\langle f,\omega_s\rangle+\frac{1}{2}\theta^2\langle f,\omega_s\rangle^2\right)\mathbb{N}_x(d\omega).$$
Then, by \eqref{N2}, we get
\begin{eqnarray}
&&\kappa(\theta_1,\theta_2,s,t)\nonumber\\
&=&\mathbb{P}_{\mu}\left(\exp\left\{i\theta_1e^{\lambda_1 t}\langle \phi_1,X_t\rangle\right.\right.\nonumber\\
&&\left.\left.+\int_{E}
\int_{\mathbb{D}}\left(-\frac{1}{2}e^{\lambda_1(t+s)}\theta^2_2\langle f, \omega_s\rangle^2
\right)\mathbb{N}_x(d\omega)X_t(dx)+\langle R_s(e^{\lambda_1(t+s)/2}\theta_2,\cdot),X_t\rangle\right\}\right)\nonumber\\
&=&\mathbb{P}_{\mu}\left(\exp\left\{i\theta_1e^{\lambda_1 t}\langle \phi_1,X_t\rangle
-\frac{1}{2}\theta^2_2e^{\lambda_1 t}\langle V_s, X_t\rangle
+\langle R_s(e^{\lambda_1(t+s)/2}\theta_2,\cdot),X_t\rangle\right\}\right),
\end{eqnarray}
where $V_s(x):=e^{\lambda_1 s}{\V}{\rm ar}_{\delta_x}\langle f, X_s\rangle.$
By \eqref{3.20}, we have
\begin{eqnarray}
\left|R_s(e^{\lambda_1(t+s)/2}\theta_2,x)\right|
&\leq& \theta_2^2e^{\lambda_1 (t+s)}\mathbb{N}_{x}
\left(\langle f,\omega_s\rangle^2\left(\frac{e^{\lambda_1(t+s)/2}\theta_2\langle f,\omega_s\rangle}{6}\wedge 1\right)\right)\nonumber\\
&=& \theta_2^2e^{\lambda_1 t}\mathbb{N}_{x}\left(Y_s^2\left(\frac{\theta_2e^{\lambda_1 t/2}Y_s}{6}\wedge 1\right)\right),
\end{eqnarray}
where $Y_s:=e^{\lambda_1 s/2}\langle f,\omega_s\rangle.$
Let
$$
h(x,s,t):=\mathbb{N}_{x}\left(Y_s^2\left(\frac{\theta_2e^{\lambda_1 t/2}Y_s}{6}\wedge 1\right)\right).
$$
We note that $h(x,s,t)\downarrow0$ as  $t\uparrow\infty$
and by \eqref{2.10}, we get
$$
h(x,s,t)\le \mathbb{N}_{x}(Y_s^2)=
e^{\lambda_1 s}\mathbb{V}{\rm ar}_{\delta_x}(\langle f,X_s\rangle)
\lesssim a_{t_0}(x)^{1/2}\in L^2(E,m).
$$
Thus, by \eqref{1.25}, we have, for any $u<t$,
$$\limsup_{t\to\infty}e^{\lambda_1t}T_t(h(\cdot,s,t))
\le \limsup_{t\to\infty}e^{\lambda_1t}T_t(h(\cdot,s,u))
=\langle h(\cdot,s,u),\phi_1\rangle_m\phi_1(x).$$
Letting $u\to\infty$, we get
$\lim_{t\to\infty}e^{\lambda_1t}T_t(h(\cdot,s,t))=0$. Therefore we have
$$
\mathbb{P}_{\mu}\left|\langle R_s(e^{\lambda_1(t+s)/2}\theta_2,\cdot),X_t\rangle\right|\le \theta^2_2e^{\lambda_1t}T_t(h(\cdot,s,t))\to 0,\quad\mbox{as } t\to\infty,
$$
which implies
$$\lim_{t\to\infty}\langle R_s(e^{\lambda_1(t+s)/2}\theta_2,\cdot),X_t\rangle=0, \quad \mbox{in probability.}$$
Furthermore, by Remark \ref{rem:large} and the fact $V_s(x)\lesssim a_{t_0}(x)^{1/2}\in L^2(E,m)\cap L^4(E,m)$, we have
\begin{equation*}
\lim_{t\to\infty}e^{\lambda_1 t}\langle V_s, X_t\rangle
=\sigma_f^2(s)W_\infty,\quad \mbox{in probability},
\end{equation*}
where
$\sigma_f^2(s):= \langle V_s,~\phi_1\rangle_m$.
Hence by the dominated convergence theorem, we get
\begin{equation}\label{10.4}
 \lim_{t\to\infty}\kappa(\theta_1,\theta_2,s,t)= \mathbb{P}_{\mu}\left(\exp\left\{i\theta_1W_\infty\right\}
  \exp\left\{-\frac{1}{2}\theta_2^2\sigma_f^2(s) W_\infty\right\}\right),
\end{equation}
which implies our claim \eqref{10.5}.

Since $e^{\lambda_1 (t+s)}\langle \phi_1,X_{t+s}\rangle-e^{\lambda_1 t}\langle \phi_1,X_t\rangle \to 0$ in probability, as $t\to\infty$, we easily get that under $\mathbb{P}_\mu$,
\begin{eqnarray*}
U_3(s,t):=\left(e^{\lambda_1 (t+s)}\langle \phi_1,X_{t+s}\rangle,e^{\lambda_1(t+s)/2}(\langle f,X_{t+s}\rangle-\langle T_sf,X_{t}\rangle)\right)
\stackrel{d}{\to}(W_\infty, \sqrt{W_\infty}G_1(s)),
\end{eqnarray*}
 as  $t\to\infty$.
By \eqref{2.2}, we have $\lim_{s\to\infty}V_s(x)= \sigma^2_f\phi_1(x)$,
thus $\lim_{s\to\infty}\sigma_f^2(s)=\sigma^2_f$. So
\begin{equation}\label{10.7}
\lim_{s\to\infty}d(G_1(s),G_1(f))=0.
\end{equation}
Let $\mathcal{D}(s+t)$ and $\widetilde{\mathcal{D}}(s,t)$ be the distributions of $U_1(s+t)$ and $U_3(s,t)$ respectively,
and let $\mathcal{D}(s)$ and $\mathcal{D}$ be the distributions of $(W_\infty, \sqrt{W_\infty}G_1(s))$
and $(W_\infty, \sqrt{W_\infty}G_1(f))$ respectively.
Then, using \eqref{5.20}, we have
\begin{eqnarray}\label{10.11}
  \limsup_{t\to\infty}d(\mathcal{D}(s+t),\mathcal{D})&\leq&
  \limsup_{t\to\infty}[d(\mathcal{D}(s+t),\widetilde{\mathcal{D}}(s,t))+d(\widetilde{\mathcal{D}}(s,t),\mathcal{D}(s))
  +d(\mathcal{D}(s),\mathcal{D})]\nonumber\\
 &\leq &\limsup_{t\to\infty}(\mathbb{P}_{\mu}(e^{\lambda_1(t+s)/2}\langle T_s f,X_{t}\rangle )^2)^{1/2}+0+d(\mathcal{D}(s),\mathcal{D}).
\end{eqnarray}
Using this and the definition of $\limsup_{t\to\infty}$, we easily get that
$$
\limsup_{t\to\infty}d (\mathcal{D}(t),\mathcal{D})=\limsup_{t\to\infty}d(\mathcal{D}(s+t),\mathcal{D})
\le \limsup_{t\to\infty}(\mathbb{P}_{\mu}(e^{\lambda_1(t+s)/2}\langle T_s f,X_{t}\rangle)^2)^{1/2}+d(\mathcal{D}(s),\mathcal{D}).
$$
Letting $s\to\infty$, we get
$$ \limsup_{t\to\infty}d(\mathcal{D}(t),\mathcal{D})
\le\limsup_{s\to\infty}\limsup_{t\to\infty}(\mathbb{P}_{\mu}(e^{\lambda_1(t+s)/2}\langle T_s f,X_{t}\rangle)^2)^{1/2}.$$
Therefore, we are left to prove that
\begin{equation}\label{limsuplimsup}
\limsup_{s\to\infty}\limsup_{t\to\infty}e^{\lambda_1(t+s)}\mathbb{P}_{\mu}(\langle T_s f,X_{t}\rangle)^2=0.
\end{equation}

By \eqref{1.13} and \eqref{1.36}, we have for any $x\in E$,
\begin{eqnarray*}
  &&e^{\lambda_1(t+s)}\mathbb{V}{\rm ar}_{\delta_x}\langle T_s f,X_{t}\rangle
  = e^{\lambda_1(s+t)}\int_0^t T_{t-u}[A(T_{s+u}f)^2](x)\,du\\
  &=&e^{(\lambda_1-2\lambda_{\gamma(f)})s}\int_0^t e^{(\lambda_1-2\lambda_{\gamma(f)})u}e^{\lambda_1(t-u)}T_{t-u}[A(e^{\lambda_{\gamma(f)}(s+u)}T_{s+u}f)^2](x)\,du\\
  &\lesssim& e^{(\lambda_1-2\lambda_{\gamma(f)})s}
  \left(\int_0^t e^{(\lambda_1-2\lambda_{\gamma(f)})u}e^{\lambda_1(t-u)}T_{t-u}[a_{2t_0}](x)\,du\right)
\end{eqnarray*}
and
\begin{eqnarray*}
  &&\int_0^t e^{(\lambda_1-2\lambda_{\gamma(f)})u}e^{\lambda_1(t-u)}T_{t-u}(a_{2t_0})(x)\,du\\
  &=&\left(\int_0^{t-t_0}+\int_{t-t_0}^t\right)
  e^{(\lambda_1-2\lambda_{\gamma(f)})u}e^{\lambda_1(t-u)}T_{t-u}(a_{2t_0})(x)\,du\\
  &\lesssim& \int_0^{t-t_0}e^{(\lambda_1-2\lambda_{\gamma(f)})u}\,du a_{t_0}(x)^{1/2}
  +\int_0^{t_0}e^{(\lambda_1-2\lambda_{\gamma(f)})(t-u)}e^{\lambda_1u}T_{u}(a_{2t_0})(x)\,du\\
  &\lesssim& a_{t_0}(x)^{1/2}+\int_0^{t_0}T_{u}(a_{2t_0})(x)\,du\lesssim a_{t_0}(x)^{1/2}.
\end{eqnarray*}
The last inequality follows from \eqref{1.19}.
Thus,
\begin{eqnarray}\label{2.15}
  \limsup_{t\to\infty}e^{\lambda_1(t+s)}\mathbb{V}{\rm ar}_{\mu}\langle T_s f,X_{t}\rangle
  &=&\limsup_{t\to\infty}e^{\lambda_1(t+s)}\langle\mathbb{V}{\rm ar}_{\delta_{\cdot}}\langle T_s f,X_{t}\rangle,\mu\rangle\nonumber\\
  &\lesssim &e^{(\lambda_1-2\lambda_{\gamma(f)})s}\langle a_{t_0}(x)^{1/2},\mu\rangle.
\end{eqnarray}
By \eqref{limit-mean}, we get
\begin{equation}\label{2.14}
  \lim_{t\to\infty}e^{\lambda_1(t+s)/2}\P_{\mu}\langle T_s f,X_{t}\rangle=\lim_{t\to\infty}e^{\lambda_1(t+s)/2}\langle T_{(t+s)}f, \mu\rangle=0.
\end{equation}
Now \eqref{limsuplimsup} follows easily from \eqref{2.15} and \eqref{2.14}. The proof is now complete.
\hfill$\Box$

\bigskip

\begin{lemma}\label{lem:5.5}
Assume that $f\in \C_s$ and $h\in\C_c$.
Define
$$
  Y_1(t):=t^{-1/2}e^{\lambda_1t/2}\langle h,X_t\rangle, \quad Y_2(t):=e^{\lambda_1t/2}\langle f,X_t\rangle, \quad t>0,
$$
and
$$Y_t:=Y_1(t)+Y_2(t).
$$
Then for any $c>0$, $\delta>0$ and $x\in E$, we have
\begin{equation}\label{4.5}
 \lim_{t\to\infty}\mathbb{P}_{\delta_x}\left(|Y_t|^2;|Y_t|>ce^{\delta t }\right)=0.
\end{equation}
\end{lemma}
\textbf{Proof:}
For any $\epsilon>0$ and $\eta>0$, we have
\begin{eqnarray*}
  \mathbb{P}_{\delta_x}\left(|Y_t|^2;|Y_t|>ce^{\delta t }\right)
  &\le& 2\mathbb{P}_{\delta_x}\left(|Y_1(t)|^2;|Y_t|>ce^{\delta t }\right)
  +2\mathbb{P}_{\delta_x}\left(|Y_2(t)|^2;|Y_t|>ce^{\delta t }\right) \\
   &\le& 2\mathbb{P}_{\delta_x}\left(|Y_1(t)|^2;|Y_1(t)|>\epsilon e^{\delta t }\right)
   +2 \epsilon^2 e^{2\delta t}\mathbb{P}_{\delta_x}\left(|Y_t|>c e^{\delta t }\right)\\
   && +2\mathbb{P}_{\delta_x}\left(|Y_2(t)|^2;|Y_2(t)|^2>\eta\right)+2\eta\mathbb{P}_{\delta_x}\left(|Y_t|>ce^{\delta t }\right)\\
   &=:& J_1(t,\epsilon)+J_2(t,\epsilon)+J_3(t,\eta)+J_4(t,\eta).
\end{eqnarray*}
Repeating the proof of \cite[Lemma 3.2]{RSZ2} (with the $S_tf$ there replaced by $Y_1(t)$), we can get
\begin{equation}\label{J1}
  \lim_{t\to\infty}J_1(t,\epsilon)=2\lim_{t\to\infty}\mathbb{P}_{\delta_x}\left(|Y_1(t)|^2;|Y_1(t)|>\epsilon e^{\delta t }\right)=0.
\end{equation}
By \eqref{limit-mean} and \eqref{2.2}, we easily get
\begin{equation}\label{2.16}
  \lim_{t\to\infty}\mathbb{P}_{\delta_x}(|Y_2(t)|^2)=\sigma^2_f\phi_1(x).
\end{equation}
By \eqref{1.49} and the fact $\mathbb{P}_{\delta_x}(Y_1(t))=t^{-1/2}h(x)$, we get
$$
\lim_{t\to\infty}\mathbb{P}_{\delta_x}(|Y_1(t)|^2)=\lim_{t\to\infty}\left(\mathbb{V}{\rm ar}_{\delta_x}(Y_1(t))+t^{-1}h^2(x)\right)=\rho^2_h\phi_1(x).
$$
Thus,
\begin{equation}\label{2.17}
  \limsup_{t\to\infty}\mathbb{P}_{\delta_x}(|Y_t|^2)\le  2\lim_{t\to\infty}\mathbb{P}_{\delta_x}(|Y_1(t)|^2+|Y_2(t)|^2)=2(\sigma^2_f+\rho^2_h)\phi_1(x).
\end{equation}
Thus by Chebyshev's inequality, we have
\begin{eqnarray}\label{J2}
  \lim_{\epsilon\to0}\limsup_{t\to\infty}J_2(t,\epsilon) &\le& 2\lim_{\epsilon\to0}\epsilon^2c^{-2}\limsup_{t\to\infty}\mathbb{P}_{\delta_x}(|Y_t|^2)=0.
\end{eqnarray}

For $J_3(t,\eta)$, by Lemma \ref{lem:small}, $Y_2(t)\stackrel{d}{\to}G_1(f)\sqrt{W_\infty}$.
Let $\Psi_{\eta}(r)=r$ on $[0,\eta-1]$,
$\Psi_\eta(r)=0$ on $[\eta,\infty]$, and let $\Psi_\eta$ be linear on $[\eta-1,\eta]$.
Then, by \eqref{2.16},
\begin{eqnarray*}
  \limsup_{t\to\infty}\mathbb{P}_{\delta_x}\left(|Y_2(t)|^2;|Y_2(t)|^2>\eta\right) &=& \limsup_{t\to\infty}\left(\mathbb{P}_{\delta_x}\left(|Y_2(t)|^2\right)-\mathbb{P}_{\delta_x}
  \left(|Y_2(t)|^2;|Y_2(t)|^2\le \eta\right)\right) \\
   &\le& \limsup_{t\to\infty}\left(\mathbb{P}_{\delta_x}\left(|Y_2(t)|^2\right)-
   \mathbb{P}_{\delta_x}\left(\Psi_\eta(|Y_2(t)|^2)\right)\right)\\
   &=&\sigma^2_f\phi_1(x)-\mathbb{P}_{\delta_x}\left(\Psi_\eta(G_1(f)^2W_\infty)\right).
\end{eqnarray*}
By the monotone convergence theorem and the fact that $G_1(f)$ and $W_\infty$ are independent, we have
$$
\lim_{\eta\to\infty}\mathbb{P}_{\delta_x}\left(\Psi_\eta(G_1(f)^2W_\infty)\right)
=\mathbb{P}_{\delta_x}\left(G_1(f)^2W_\infty\right)
=\mathbb{P}_{\delta_x}\left(G_1(f)^2\right)\mathbb{P}_{\delta_x}W_\infty=\sigma^2_f\phi_1(x).
$$
Thus,
\begin{equation}\label{J3}
  \lim_{\eta\to\infty}\limsup_{t\to\infty}J_3(t,\eta)=0.
\end{equation}

By Chebyshev's inequality and \eqref{2.17},
\begin{equation}\label{J4}
  \limsup_{t\to\infty}J_4(t,\eta)\le2\eta c^{-2}\limsup_{t\to\infty}e^{-2\delta t}\mathbb{P}_{\delta_x}(|Y_t|^2)=0.
\end{equation}
Thus, \eqref{4.5} follows easily from \eqref{J1}, \eqref{J2}, \eqref{J3} and \eqref{J4}.
\hfill$\Box$

\smallskip

\begin{lemma}\label{lem:5.6}
Assume that $f\in \C_s$ and $h\in\C_c$.
Define
$$
  \widetilde{Y}_1(t)(\omega):=t^{-1/2}e^{\lambda_1t/2}\langle h,\omega_t\rangle, \quad \widetilde{Y}_2(t)(\omega):=e^{\lambda_1t/2}\langle f,\omega_t\rangle, \quad t>0,\omega\in \D,
$$
and
$$\widetilde{Y}_t:=\widetilde{Y}_1(t)+\widetilde{Y}_2(t).
$$
For any $c>0$ and $\delta>0$, we have
\begin{equation}\label{4.41}
\lim_{t\to\infty}\mathbb{N}_x\left(|\widetilde{Y}_t|^2;|\widetilde{Y}_t|>ce^{\delta t }\right)=0.
\end{equation}
\end{lemma}

\textbf{Proof:}\quad
For $t>1$,
$$
\mathbb{N}_x\left(|\widetilde{Y}_t|^2;|\widetilde{Y}_t|>ce^{\delta t }\right)
=\mathbb{N}_x\left(|\widetilde{Y}_t|^2;|\widetilde{Y}_t|>ce^{\delta t },\|\omega_{1}\|\neq 0\right).
$$
Thus, we only need to prove
$$
\lim_{t\to\infty}\widetilde{\mathbb{N}}_x\left(|\widetilde{Y}_t|^2;|\widetilde{Y}_t|>ce^{\delta t }\right)=0.
$$
For any $x\in E$,
let $N(d\omega)$ be a Poisson random measure with intensity $\mathbb{N}_x(d\omega)$ defined on the probability space $\{\tilde{\Omega},\tilde{\mathcal{F}},\mathbf{P}_{\delta_x}\}$
and
$$
\Lambda_t=\int_{\mathbb{D}}\omega_tN(d\omega).
$$
We know  that, under $\mathbf{P}_{\delta_x}$, $\{\Lambda_t,t\ge0\}$ has the same law as $\{X_t,t\ge0\}$ under $\P_{\delta_x}$.
Define
$$
\Lambda^*_t:=\int_{\widetilde{\mathbb{D}}}\omega_tN(d\omega) \quad \mbox{and }
Y_t(\Lambda^*):=t^{-1/2}e^{\lambda_1t/2}\langle h, \Lambda_t^*\rangle+e^{\lambda_1t/2}\langle f, \Lambda_t^*\rangle,
$$
where $\widetilde{\mathbb{D}}:=\{\omega\in \D: \|\omega_{1}\|\neq 0\}.$
It is clear that for $t>1$, $\Lambda_t^*=\Lambda_t$ and $Y_t(\Lambda^*)\stackrel{d}{=}Y_t$.
Since $\mathbb{N}_x(\widetilde{\mathbb{D}})<\infty$, $\Lambda_t^*$ is a compound Poisson process and can be written as
$$
\Lambda_t^*=\sum_{j=1}^K \widetilde{X}^j_{t},
$$
where $\widetilde{X}^j_{t},j=1,2,\dots$ are i.i.d.
with the same law as $\omega_t$ under $\tilde{\mathbb{N}}_x$
and $K$ is a Poisson random variable with parameter
 $\mathbb{N}_x(\widetilde{\mathbb{D}})$ which is independent of $\widetilde{X}^j_{t},j=1,2,\dots$ .
Let
$$
Y_t(\widetilde{X}^j):=t^{-1/2}e^{\lambda_1t/2}\langle h, \widetilde{X}^j_{t}\rangle+e^{\lambda_1t/2}\langle f, \widetilde{X}^j_{t}\rangle.
$$
Then, $Y_t(\widetilde{X}^j)$ is independent of $K$ and has the same law as $\widetilde{Y}_t$ under $\widetilde{\mathbb{N}}_x$.
Therefore, for $t>1$,
\begin{eqnarray*}
  \P_{\delta_x}(|Y_t|^2;|Y_t|>ce^{\delta t})&=& \mathbf{P}_{\delta_x} (|Y_t(\Lambda^*)|^2;|Y_t(\Lambda^*)|>ce^{\delta t})\\
   &\ge&  \mathbf{P}_{\delta_x} (|Y_t(\widetilde{X}^1)|^2;|Y_t(\widetilde{X}^1)|>ce^{\delta t}, K=1)\\
   &=& \mathbf{P}_{\delta_x} (K=1)\mathbf{P}_{\delta_x} (|Y_t(\widetilde{X}^1)|^2;|Y_t(\widetilde{X}^1)|>ce^{\delta t})\\
   &=&  \mathbb{N}_x(\widetilde{\mathbb{D}})e^{-\mathbb{N}_x(\widetilde{\mathbb{D}})} \widetilde{\mathbb{N}}_x(|\widetilde{Y}_t|^2;|\widetilde{Y}_t|>ce^{\delta t}).
\end{eqnarray*}
Now \eqref{4.41} follows easily from Lemma \ref{lem:5.5}.

\hfill$\Box$

\begin{lemma}\label{lem:cs}
Assume that $f\in \C_s$ and $h\in\C_c$.
Then
\begin{equation}\label{cs}
  \left(e^{\lambda_1t}\langle\phi_1,X_t\rangle, t^{-1/2}e^{\lambda_1t/2}\langle h,X_t\rangle,
  e^{\lambda_1t/2}\langle f,X_t\rangle\right)
  \stackrel{d}{\to}\left(W_\infty, \sqrt{W_\infty}G_2(h),\sqrt{W_\infty}G_1(f)\right),
\end{equation}
where $G_2(h)\sim\mathcal{N}(0,\rho^2_h)$ and $G_1(f)\sim\mathcal{N}(0,\sigma^2_f)$.
Moreover, $W_\infty$, $G_2(h)$ and $G_1(f)$ are independent.

\end{lemma}
\textbf{Proof:}\quad
In the proof, we always assume $t>3t_0$.
We define an ${\mathbb R}^3$-valued random variable by
\begin{equation*}
  U_1(t):=\left(e^{\lambda_1t}\langle\phi_1,X_t\rangle, t^{-1/2}e^{\lambda_1t/2}\langle h,X_t\rangle,
  e^{\lambda_1t/2}\langle f,X_t\rangle\right).
\end{equation*}
Let $n>2$ and write
\begin{equation*}
  U_1(nt)=\left(e^{\lambda_1nt}\langle\phi_1,X_{nt}\rangle,\, (nt)^{-1/2}e^{\lambda_1nt/2}\langle h,X_{nt}\rangle,
  e^{\lambda_1nt/2}\langle f,X_{nt}\rangle\right).
\end{equation*}
To consider the limit of  $U_1(t)$ as $t\to\infty$, it is equivalent to
consider  the limit of  $U_1(nt)$ for any $n>2$. The main idea is as
follows. For $t>t_0, n>2$,
\begin{equation}\label{decom-U1(nt)}
  \begin{array}{rll}U_1(nt)&=&\displaystyle\left(e^{\lambda_1nt}\langle\phi_1,X_{nt}\rangle,\,
\frac{e^{\lambda_1nt/2}(\langle h,X_{nt}\rangle-\langle T_{(n-1)t}h,X_t\rangle)}{((n)t)^{1/2}},\,
e^{\lambda_1nt/2}(\langle f,X_{nt}\rangle-\langle T_{(n-1)t}f,X_t\rangle)\right)\\
&&\displaystyle+\left(0,\,
(nt)^{-1/2}e^{\lambda_1nt/2}\langle T_{(n-1)t}h,X_t\rangle,\,
e^{\lambda_1nt/2}\langle T_{(n-1)t}f,X_t\rangle\right).
\end{array}
\end{equation}
The double limit, first as $t\to\infty$ and then $n\to \infty$, of the first term of the right side of \eqref{decom-U1(nt)}
is equal to the double limit, first as $t\to\infty$ and then $n\to \infty$, of
another ${\mathbb R}^2$-valued random variable $U_2(n,t)$
where
\begin{eqnarray*}
&&U_2(n,t)\\
&:=&\left(e^{\lambda_1t}\langle\phi_1,X_{t}\rangle,\,
\frac{e^{\lambda_1nt/2}(\langle h,X_{nt}\rangle-\langle T_{(n-1)t}h,X_t\rangle)}{((n-1)t)^{1/2}},\,
e^{\lambda_1nt/2}(\langle f,X_{nt}\rangle-\langle T_{(n-1)t}f,X_t\rangle)\right).
\end{eqnarray*}
We will prove that the second term on the right hand side of \eqref{decom-U1(nt)} has no contribution to
the double limit, first as $t\to\infty$ and then $n\to \infty$, of the left
hand side of \eqref{decom-U1(nt)}.

We claim that
\begin{equation}\label{9.5}
 U_2(n,t)\stackrel{d}{\to}\left(W_\infty, \sqrt{W_\infty}G_2(h),\sqrt{W_\infty}G_1(f) \right), \quad \mbox{ as } t\to\infty.
\end{equation}
Denote the characteristic function of $U_2(n,t)$ under $\P_\mu$ by
$\kappa_2(\theta_1,\theta_2,\theta_3,n,t)$.
Define
$$
  Y_1(t,\theta_2):=\theta_2t^{-1/2}e^{\lambda_1t/2}\langle h,X_t\rangle, \quad Y_2(t,\theta_3):=\theta_3e^{\lambda_1t/2}\langle f,X_t\rangle, \quad t>0,
$$
and
$$Y_t(\theta_2,\theta_3)=Y_1(t,\theta_2)+Y_2(t,\theta_3).$$
We define the corresponding random variables on $\D$ as $\tilde{Y}_1(t,\theta_2),\tilde{Y}_2(t,\theta_3)$ and $\tilde{Y}_t(\theta_2,\theta_3)$.
Using an argument similar to that leading to \eqref{10.9}, we get
\begin{eqnarray*}
  \kappa_2(\theta_1,\theta_2,\theta_3,n,t)
  &=&\mathbb{P}_{\mu}\left(\exp\left\{i\theta_1e^{\lambda_1 t}\langle \phi_1,X_t\rangle
  +\int_{E}\int_{\mathbb{D}}
     \left(\exp\left\{ ie^{\lambda_1 t/2}\tilde{Y}_{(n-1)t}(\theta_2,\theta_3)(\omega)\right\}\right.\right.\right.\nonumber\\
  &&\left.\left.\left.
  -1-ie^{\lambda_1 t/2}\tilde{Y}_{(n-1)t}(\theta_2,\theta_3)(\omega)\right)\mathbb{N}_x(d\omega)X_t(dx)\right\}\right).
\end{eqnarray*}
Define
\begin{equation*}
   R_t'(x,\theta):=\int_{\mathbb{D}}\left(\exp\{ i\theta \tilde{Y}_{t}(\theta_2,\theta_3)(\omega)\}-1-i\theta \tilde{Y}_{t}(\theta_2,\theta_3)(\omega)
   +\frac{1}{2}\theta^2 (\tilde{Y}_{t}(\theta_2,\theta_3)(\omega))^2\right)\mathbb{N}_x(d\omega)
\end{equation*}
and
$$
J(n,t,x):=\int_{\mathbb{D}}\left(\exp\{ ie^{\lambda_1 t/2}\tilde{Y}_{(n-1)t}(\theta_2,\theta_3)(\omega)\}-1-ie^{\lambda_1 t/2}\tilde{Y}_{(n-1)t}(\theta_2,\theta_3)(\omega)\right)\mathbb{N}_x(d\omega).
$$
Then
\begin{eqnarray*}
  J(n,t,x)
  =-\frac{1}{2}e^{\lambda_1 t}\mathbb{N}_x(\tilde{Y}_{(n-1)t}(\theta_2,\theta_3))^2+R'_{(n-1)t}(x,e^{\lambda_1 t/2}),
\end{eqnarray*}
and
\begin{equation*}
  \kappa_2(\theta_1,\theta_2,\theta_3,n,t)
  =\mathbb{P}_{\mu}\left(\exp\left\{i\theta_1e^{\lambda_1 t}\langle \phi_1,X_t\rangle+\langle J(n,t,\cdot),X_t\rangle\right\}\right).
\end{equation*}
Let $V^n_t(x):=\mathbb{N}_x(\tilde{Y}_{(n-1)t}(\theta_2,\theta_3))^2$. Then
\begin{eqnarray*}
   \langle J(n,t,\cdot),X_t\rangle&=&-\frac{1}{2}e^{\lambda_1 t}\langle V^n_t,X_t\rangle
   +\langle R'_{(n-1)t}(\cdot,e^{\lambda_1 t/2}),X_t\rangle\\
   &:=&J_1(n,t)+J_2(n,t).
\end{eqnarray*}

We first consider $J_1(n,t)$.
By \eqref{N2},
\begin{eqnarray*}
  &&V^n_t(x)=\mathbb{V}{\rm ar}_{\delta_x}(Y_{(n-1)t}(\theta_2,\theta_3))\\
  &=&\mathbb{V}{\rm ar}_{\delta_x}(Y_1((n-1)t,\theta_2))+\mathbb{V}{\rm ar}_{\delta_x}(Y_2((n-1)t,\theta_3))+ \mathbb{C}{\rm ov}_{\delta_x}(Y_1((n-1)t,\theta_2),Y_2((n-1)t,\theta_3)).
\end{eqnarray*}
So by \eqref{1.49}, \eqref{small:new} and \eqref{cov:sc}, we have, for $t>3t_0$,
\begin{eqnarray}\label{8.10}
  &&\left|V^n_t(x)-(\theta_2^2\rho_h^2+\theta_3^2\sigma_f^2)\phi_1(x)\right| \nonumber\\
  &\le& \left|\mathbb{V}{\rm ar}_{\delta_x}(Y_1((n-1)t,\theta_2))-\theta_2^2\rho_h^2\phi_1(x)\right|
  +\left|\mathbb{V}{\rm ar}_{\delta_x}(Y_2((n-1)t,\theta_3))-\theta_3^2\sigma_f^2\phi_1(x)\right|\nonumber\\
  &&+\left|\mathbb{C}{\rm ov}_{\delta_x}(Y_1((n-1)t,\theta_2),Y_2((n-1)t,\theta_3))\right|\nonumber\\
  &\lesssim&\left(e^{(\lambda_1-2\lambda_{\gamma(f)})(n-1)t}+e^{(\lambda_1-\lambda_2)(n-1)t}+((n-1)t)^{-1/2}+((n-1)t)^{-1}\right)a_{t_0}(x)^{1/2}.
\end{eqnarray}
Thus, we have that as $t\to\infty$,
\begin{eqnarray*}
   &&e^{\lambda_1 t}\langle \left|V^n_t(x)-(\theta_2^2\rho_h^2+\theta_3^2\sigma_f^2)\phi_1(x)\right|, X_t\rangle\\
&\lesssim&  \left(e^{(\lambda_1-2\lambda_{\gamma(f)})(n-1)t}+e^{(\lambda_1-\lambda_2)(n-1)t}+((n-1)t)^{-1/2}+((n-1)t)^{-1}\right)
   e^{\lambda_1 t}\langle (a_{t_0})^{1/2}, X_t\rangle\to 0,
\end{eqnarray*}
in probability.
It follows that
\begin{equation}\label{8.15}
  \lim_{t\to\infty}J_1(n,t)=\lim_{t\to\infty}
  -\frac{1}{2}e^{\lambda_1 t}(\theta_2^2\rho_h^2+\theta_3^2\sigma_f^2)\langle \phi_1, X_t\rangle
  =-\frac{1}{2}(\theta_2^2\rho_h^2+\theta_3^2\sigma_f^2)W_\infty \quad \mbox{in probability}.
\end{equation}

For $J_{2}(n,t)$, by \eqref{3.20}, we have, for any $\epsilon>0$,
\begin{eqnarray*}
|R'_{(n-1)t}(x,e^{\lambda_1t/2})|&\leq& \frac{1}{6}e^{\frac{3}{2}\lambda_1 t}
 \mathbb{N}_x\left(|\tilde{Y}_{(n-1)t}(\theta_2,\theta_3)|^3;
|\tilde{Y}_{(n-1)t}(\theta_2,\theta_3)|<\epsilon e^{-\lambda_1 t/2}\right)\\
&&+e^{\lambda_1 t}\mathbb{N}_x
\left(|\tilde{Y}_{(n-1)t}(\theta_2,\theta_3)|^2; |\tilde{Y}_{(n-1)t}(\theta_2,\theta_3)|\geq\epsilon e^{-\lambda_1 t/2}\right)\\
&\leq&  \frac{\epsilon}{6}e^{\lambda_1 t}
\mathbb{N}_x\left(|\tilde{Y}_{(n-1)t}(\theta_2,\theta_3)|^2\right)\\
&&+e^{\lambda_1 t}\mathbb{N}_x
\left(|\tilde{Y}_{(n-1)t}(\theta_2,\theta_3)|^2; |\tilde{Y}_{(n-1)t}(\theta_2,\theta_3)|\geq\epsilon e^{-\lambda_1 t/2}\right)\\
&=&\frac{\epsilon}{6}e^{\lambda_1 t}V^n_t(x)+e^{\lambda_1 t}F_t^n(x),
\end{eqnarray*}
where
$F^n_t(x)=\mathbb{N}_x\left(|\tilde{Y}_{(n-1)t}(\theta_2,\theta_3)|^2; |\tilde{Y}_{(n-1)t}(\theta_2,\theta_3)|\geq\epsilon e^{-\lambda_1 t/2}\right)$.
Note that
\begin{equation}\label{mean-F}
e^{\lambda_1 t} \P_{\mu}\langle F_t^n(x), X_t\rangle=e^{\lambda_1t}\langle T_t(F^n_t),\mu\rangle.
\end{equation}
It follows from Lemma \ref{lem:5.6} that $\lim_{t\to\infty}F^n_t(x)= 0$.
By \eqref{8.10}, we also have
$$
F^n_t(x)\leq V^n_t(x)\lesssim a_{t_0}(x)^{1/2},
$$
which implies that
$\|F^n_t\|_2\to 0$ as $t\to\infty$.
By Lemma \eqref{lem:rsnew},
$$
\lim_{t\to\infty}e^{\lambda_1t}T_t(F^n_t)(x)=0.
$$
Note that,
by \eqref{1.36},
$e^{\lambda_1t}T_t(F^n_t)\lesssim e^{\lambda_1t}T_t(a^{1/2}_{t_0})\lesssim a^{1/2}_{t_0}$.
Since $\mu$ has compact support and $a_{t_0}$ is continuous, we have $\langle a_{t_0},\mu\rangle<\infty$.
By \eqref{mean-F} and the dominated convergence theorem, we obtain
$\lim_{t\to\infty}e^{\lambda_1 t} \P_{\mu}\langle F_t^n(x), X_t\rangle=0$, which implies that
$e^{\lambda_1 t}\langle F_t^n(x), X_t\rangle\to 0$ in probability.
Furthermore, by \eqref{8.15}, we have that as $t\to\infty$,
$$
\frac{\epsilon}{6}e^{\lambda_1 t}\langle V^n_t,
X_t\rangle\to\frac{\epsilon}{6} (\theta_2^2\rho_h^2+\theta_3^2\sigma_f^2)W_\infty\quad \mbox{in probability}.
$$
Thus, letting $\epsilon\to 0$, we get that as $t\to\infty$,
\begin{equation}\label{4.48}
 J_2(n,t)\to0 \quad \mbox{in probability}.
\end{equation}
Thus, when $t\to\infty$,
\begin{equation}\label{9.8}
  \exp\left\{\langle J(n,t,\cdot),X_t\rangle\right\}
  \to\exp\left\{-\frac{1}{2}(\theta_2^2\rho_h^2+\theta_3^2\sigma_f^2) W_\infty\right\}
\end{equation}
 in probability. Since the real part of $J(n,t,x)$ is less than 0,
 we have
$$
|\exp\left\{\langle J(n,t,\cdot),X_t\rangle\right\}|\leq 1.
$$
So by the dominated convergence theorem, we get that
\begin{equation}\label{9.10}
\lim_{t\to\infty}\kappa_2(\theta_1,\theta_2,\theta_3,n,t)=\P_{\mu}\left[\exp\left\{i\theta_1W_\infty\right\}
  \exp\left\{-\frac{1}{2}(\theta_2^2\rho_h^2+\theta_3^2\sigma_f^2)  W_\infty\right\}\right],
 \end{equation}
which implies our claim \eqref{9.5}.

By \eqref{9.5}
and the fact $e^{\lambda_1 nt}\langle \phi_1,X_{nt}\rangle-e^{\lambda_1 t}\langle \phi_1,X_{t}\rangle\to 0$, in probability,
as $t\to\infty$ , we easily get
\begin{eqnarray*}
&&U_3(n,t)\\
&&:=\left(e^{\lambda_1 nt}\langle \phi_1,X_{nt}\rangle,
\frac{e^{\lambda_1nt/2}(\langle h,X_{nt}\rangle-\langle T_{(n-1)t}h,X_t\rangle)}{(nt)^{1/2}},
e^{\lambda_1nt/2}(\langle f,X_{nt}\rangle-\langle T_{(n-1)t}f,X_t\rangle)\right)\\
&&\stackrel{d}{\to}\left(W_\infty, \sqrt{\frac{n-1}{n}}\sqrt{W_\infty}G_2(h),\sqrt{W_\infty}G_1(f)\right).
\end{eqnarray*}
Using \eqref{1.49} and the fact
${\P}_{\mu}\langle h,X_t\rangle=\langle T_th,\mu\rangle= e^{-\lambda_1 t/2}\langle h,\mu\rangle$,
we can get
\begin{eqnarray}\label{9.2}
  (nt)^{-1}e^{\lambda_1 nt}{\P}_{\mu}(\langle T_{(n-1)t}h,X_t\rangle)^2
 &=&(nt)^{-1}e^{\lambda_1 t}\mathbb{V}ar_{\mu}\langle h,X_t\rangle
  +(nt)^{-1}e^{\lambda_1 t}({\P}_{\mu}\langle h,X_t\rangle)^2\nonumber\\
  &\lesssim& n^{-1}(1+t^{-1}).
\end{eqnarray}
Using \eqref{2.15} with $s=(n-1)t$, and then letting $t\to\infty$, by \eqref{limit-mean} we get
\begin{equation}\label{9.3}
  e^{\lambda_1 nt}{\P}_{\mu}\langle T_{(n-1)t}f,X_t\rangle)^2\lesssim
e^{(\lambda_1-2\lambda_{\gamma(f)})(n-1)t}\langle a_{t_0}(x)^{1/2},\mu\rangle+e^{\lambda_1 nt}\langle T_{nt}f,\mu\rangle^2\to0.
\end{equation}
Let $\mathcal{D}(nt)$ and $\widetilde{\mathcal{D}}^n(t)$ be the distributions of
$U_1(nt)$ and $U_3(n,t)$ respectively,
and let $\mathcal{D}^n$ and $\mathcal{D}$ be the distributions of
$\left(W_\infty, \sqrt{\frac{n-1}{n}}\sqrt{W_\infty}G_2(h),\sqrt{W_\infty}G_1(f)\right)$ and $\left(W_\infty, \sqrt{W_\infty}G_2(h),\sqrt{W_\infty}G_1(f)\right)$ respectively.
Then, using \eqref{5.20}, we have
\begin{eqnarray}\label{4.12}
 &&\limsup_{t\to\infty}d(\mathcal{D}(nt),\mathcal{D})\leq
  \limsup_{t\to\infty}[d(\mathcal{D}(nt),\widetilde{\mathcal{D}}^n(t))
  +d(\widetilde{\mathcal{D}}^n(t),\mathcal{D}^n)+d(\mathcal{D}^n,\mathcal{D})]\nonumber\\
 &\leq &\limsup_{t\to\infty}\left((nt)^{-1}e^{\lambda_1 nt}{\P}_{\mu}\langle T_{(n-1)t}h,X_t\rangle^2+e^{\lambda_1 nt}{\P}_{\mu}\langle T_{(n-1)t}f,X_t\rangle^2\right)^{1/2}
 +0+d(\mathcal{D}^n,\mathcal{D}).\nonumber\\
\end{eqnarray}
Using  the definition of
 $\limsup_{t\to\infty}$, \eqref{9.2} and \eqref{9.3}, we easily get that
$$
\limsup_{t\to\infty}d(\mathcal{D}(t),\mathcal{D})=
\limsup_{t\to\infty}d(\mathcal{D}(nt),\mathcal{D})
\le c/\sqrt{n}+d(\mathcal{D}^n,\mathcal{D}),
$$
where $c$ is a constant.
Letting $n\to\infty$, we get $ \limsup_{t\to\infty}d(\mathcal{D}(t),\mathcal{D})=0$.
The proof is now complete.

\hfill$\Box$

Recall that
$$g(x)=\sum_{k: 2\lambda_k<\lambda_1}\sum_{j=1}^{n_k}b_j^k\phi_j^{(k)}(x)\quad \mbox{and }\quad I_ug(x)=\sum_{k: 2\lambda_k<\lambda_1}\sum_{j=1}^{n_k}e^{\lambda_ku}b_j^k\phi_j^{(k)}(x).$$
Note that the sum over $k$ is a sum over a finite number of elements.
Define
$$
H_\infty(\omega):=\sum_{k: 2\lambda_k<\lambda_1}\sum_{j=1}^{n_k}b_j^k\widetilde{H}_\infty^{k,j}(\omega),\quad \omega\in \D.
$$
By Lemma \ref{lem:2.5}, we have, as $u\to\infty$
$$
\langle I_ug,\omega_u\rangle\to H_\infty,\quad \mathbb{N}_x\mbox{-a.e.},
\quad \mbox{in }L^1(\mathbb{N}_x)\quad \mbox{and in } L^2(\mathbb{N}_x).
$$
Since $\mathbb{N}_x\langle I_ug,\omega_u\rangle=\mathbb{P}_{\delta_x}\langle I_ug,X_u\rangle=g(x)$,
we get
\begin{equation}\label{L1H}
  \mathbb{N}_x(H_\infty)=g(x).
\end{equation}
By \eqref{N2} and \eqref{1.13},
we have
\begin{eqnarray}\label{var:Iu}
  \mathbb{N}_x\langle I_ug,\omega_u\rangle^2=\mathbb{V}{\rm ar}_{\delta_x}\langle I_ug,X_u\rangle
  =\int_0^uT_s\left[ A\left(\sum_{k:2\lambda_k<\lambda_1}\sum_{j=1}^{n_k}e^{\lambda_ks}b_j^k\phi_j^k\right)^2\right](x)\,ds,
\end{eqnarray}
which implies
\begin{equation}\label{L2H}
  \mathbb{N}_x(H_\infty)^2=\int_0^\infty
    T_s\left[ A\left(\sum_{k:2\lambda_k<\lambda_1}\sum_{j=1}^{n_k}e^{\lambda_ks}b_j^k\phi_j^k\right)^2\right](x)\,ds.
\end{equation}
By \eqref{1.37}, we have that for any $x\in E$,
$$
\sum_{k: 2\lambda_k<\lambda_1}\sum_{j=1}^{n_k}e^{\lambda_ks}|b_j^k||\phi_j^k(x)|\lesssim e^{\lambda_K s}a_{2t_0}(x)^{1/2},
$$
 where $K=\sup\{k:2\lambda_k<\lambda_1\}$.
So by \eqref{L2H}, \eqref{1.36} and \eqref{1.19}, we have that for any $x\in E$,

 \begin{eqnarray}\label{domi-H2}
 \mathbb{N}_x(H_\infty)^2&\lesssim& \int_0^\infty e^{(2\lambda_K-\lambda_1)s}e^{\lambda_1s}T_s(a_{2t_0})(x)\,ds\nonumber\\
    &=&\left(\int_0^{t_0}+\int_{t_0}^\infty\right) e^{(2\lambda_K-\lambda_1)s}e^{\lambda_1s}T_s(a_{2t_0})(x)\,ds\nonumber\\
    &\lesssim& \int_0^{t_0}T_s(a_{2t_0})(x)\,ds+\int_{t_0}^\infty e^{(2\lambda_K-\lambda_1)s}\,ds~a_{t_0}(x)^{1/2}\nonumber\\
    &\lesssim& a_{t_0}(x)^{1/2}\in L^2(E,m)\cap L^4(E,m).
\end{eqnarray}

Now we are ready to prove Theorem \ref{The:1.3}.

\textbf{Proof of Theorem \ref{The:1.3}:}
 Consider an $\mathbb{R}^4$-valued random variable $U_4(t)$ defined by:
 \begin{eqnarray*}\label{8.5}
   &&U_4(t)\\
   &&:= \left(e^{\lambda_1 t}\langle\phi_1,X_t\rangle,
   ~e^{\lambda_1t/2}\left(\langle g,X_t\rangle-\sum_{k: 2\lambda_k<\lambda_1}\sum_{j=1}^{n_k}e^{-\lambda_kt}b_j^kH^{k,j}_\infty\right),
   \frac{e^{\lambda_1t/2}\langle h,X_t\rangle}{t^{1/2}},
   e^{\lambda_1t/2}\langle f,X_t\rangle\right).
 \end{eqnarray*}
 To get the conclusion of Theorem \ref{The:1.3}, it suffices to show that, under $\P_{\mu}$,
\begin{equation}\label{2.5a}
   U_4(t)\stackrel{d}{\to}\left(W_\infty, \sqrt{W_\infty}G_3(g), \sqrt{W_\infty}G_2(h),\sqrt{W_\infty}G_1(f)\right),
\end{equation}
where $W_\infty$, $G_3(g)$, $G_2(h)$ and $G_1(f)$ are independent.
 Denote the characteristic function of $U_4(t)$  under  $\P_{\mu}$ by $\kappa_1(\theta_1,\theta_2,\theta_3,\theta_4,t)$.
 Then, we only need to prove
 \begin{equation}\label{8.7}
   \lim_{t\to\infty}\kappa_1(\theta_1,\theta_2,\theta_3,\theta_4,t)
   =\P_{\mu}\left(\exp\{i\theta_1W_\infty\}
   \exp\left\{-\frac{1}{2}(\theta_2^2\beta_g^2 +\theta_3^2\rho_h^2 +\theta_4^2\sigma^2_f)W_\infty\right\}\right).
 \end{equation}
 Note that, by Lemma \ref{lem:1.2}, $\sum_{k: 2\lambda_k<\lambda_1}\sum_{j=1}^{n_k}e^{-\lambda_kt}b_j^kH^{k,j}_\infty=\lim_{u\to\infty}\langle I_ug, X_{t+u}\rangle$, $\P_{\mu}$-a.s.. We have
 \begin{eqnarray}\label{8.6}
  &&\kappa_1(\theta_1,\theta_2,\theta_3,\theta_4,t)\nonumber\\
  &=&\lim_{u\to\infty}
   \P_{\mu}\left(\exp\left\{i\theta_1e^{\lambda_1 t}\langle\phi_1,X_t\rangle
   +i\theta_2e^{\lambda_1t/2}(\langle g,X_t\rangle-\langle I_ug, X_{t+u}\rangle)\right.\right.\nonumber\\
   &&\qquad\qquad\qquad\left.\left.+i\theta_3t^{-1/2}e^{\lambda_1 t/2}\langle h,X_t\rangle
   +i\theta_4e^{\lambda_1 t/2}\langle f,X_t\rangle\right\}\right)\nonumber\\
   &=&\lim_{u\to\infty}\P_{\mu}\left(\exp\left\{i\theta_1e^{\lambda_1 t}\langle\phi_1,X_t\rangle
   +i\theta_3t^{-1/2}e^{\lambda_1 t/2}\langle h,X_t\rangle
   +i\theta_4e^{\lambda_1 t/2}\langle f,X_t\rangle
   +\langle J_u(t,\cdot),X_t\rangle\right\}\right),\nonumber\\
 \end{eqnarray}
 where
 $$J_u(t,x)=\int_{\mathbb{D}}\left(\exp\left\{-i\theta_2e^{\lambda_1t/2}\langle I_ug, \omega_u\rangle\right\}
 -1+i\theta_2e^{\lambda_1t/2}\langle I_ug, \omega_u\rangle\right)\mathbb{N}_x(d\omega).$$
The last equality above follows from the Markov property of $X$, \eqref{cf} and the fact
$$\int_{\mathbb{D}}\langle I_ug,\omega_u\rangle\mathbb{N}_x(d\omega)=\P_{\delta_x}\langle I_ug,X_u\rangle=g(x).$$
We will show that
\begin{equation}\label{JutoJ}
  \lim_{u\to\infty}J_u(t,x)
  =\mathbb{N}_x\left(\exp\left\{-i\theta_2e^{\lambda_1t/2}H_\infty\right\}-1+i\theta_2e^{\lambda_1t/2}H_\infty\right)
  =: J(t,x).
\end{equation}
For $u>1$,
$|e^{- i\theta_2e^{\lambda_1t/2}\langle I_ug, \omega_u\rangle}-1|\le 2\textbf{1}_{\{\|\omega_{1}\|\neq 0\}}(\omega)$.
By Remark \ref{Rek5}, $ \mathbb{N}_x(\|\omega_{1}\|\neq 0)<\infty$.
Thus, by Lemma \ref{lem:2.5} and the dominated convergence theorem, we get
$$\lim_{u\to\infty}\int_{\mathbb{D}}\left(\exp\left\{-i\theta_2e^{\lambda_1t/2}\langle I_ug, \omega_u\rangle\right\}-1\right)\mathbb{N}_x(d\omega)
=\mathbb{N}_x\left(\exp\left\{-i\theta_2e^{\lambda_1t/2}H_\infty\right\}-1\right).$$
By \eqref{L1H}, we get
 $\mathbb{N}_xH_\infty=\mathbb{N}_x\langle I_ug,\omega_u\rangle=g(x)$. Then, \eqref{JutoJ} follows immediately .

By \eqref{3.20}, we get
$$
 \sup_{u\ge0}|J_u(t,x)|\le \frac{1}{2}\theta_2^2e^{\lambda_1t}\sup_{u\ge0}\mathbb{N}_x\langle I_ug,\omega_u\rangle^2
<\frac{1}{2}\theta_2^2e^{\lambda_1t}\mathbb{N}_xH_\infty^2<\infty.
$$
Note that, by \eqref{domi-H2},
$$\P_{\mu}\langle \mathbb{N}_{\cdot}H_\infty^2,X_t\rangle
\lesssim\P_{\mu}\langle a^{1/2}_{t_0},X_t\rangle
=\langle T_ta^{1/2}_{t_0}, \mu\rangle<\infty,$$
which implies that  $\langle \mathbb{N}_{\cdot}H_\infty^2,X_t\rangle<\infty$, $\P_\mu$-a.s.
So, by the dominated convergence theorem, we get
$$
\lim_{u\to\infty}\langle J_u(t,\cdot),X_t\rangle=\langle J(t,\cdot),X_t\rangle,\quad \P_\mu\mbox{-a.s.}
$$
Using the dominated convergence theorem again, we obtain
$$
\kappa_1(\theta_1,\theta_2,\theta_3,\theta_4,t)
=\P_{\mu}\left(\exp\left\{i\theta_1e^{\lambda_1 t}\langle\phi_1,X_t\rangle
+i\theta_3t^{-1/2}e^{\lambda_1 t/2}\langle h,X_t\rangle
   +i\theta_4e^{\lambda_1 t/2}\langle f,X_t\rangle
+\langle J(t,\cdot),X_t\rangle\right\}\right).
$$
Let
$$
R(\theta,x):=\mathbb{N}_x\left(\exp\left\{i\theta H_\infty\right\}-1-i\theta H_\infty+\frac{1}{2}\theta^2H_\infty^2\right).
$$
Thus,
$$
\langle J(t,\cdot),X_t\rangle=-\frac{1}{2}\theta_2^2e^{\lambda_1 t} \langle V,X_t\rangle+\langle R(-e^{\lambda_1 t/2}\theta_2,\cdot),X_t\rangle,
$$
where $V(x):=\mathbb{N}_x(H_\infty)^2$.
By \eqref{3.20}, we have
\begin{equation}\label{8.9}
 |R(-e^{\lambda_1 t/2}\theta_2,x)|
 \leq e^{\lambda_1 t}\theta_2^2\mathbb{N}_x\left(|H_\infty|^2\left(\frac{e^{\lambda_1 t/2}\theta_2|H_\infty|}{6}\wedge 1\right)\right),
\end{equation}
which implies that
\begin{eqnarray*}
  \P_{\mu}\left|\langle R(-e^{\lambda_1 t/2}\theta_2,\cdot),X_t\rangle\right|
  &\leq& \theta_2^2~e^{\lambda_1t}\langle~ T_t(k(\cdot,t)),~\mu~\rangle,
\end{eqnarray*}
where
$$
k(x,t):=\mathbb{N}_x\left(|H_\infty|^2\left(\frac{e^{\lambda_1t/2}\theta_2|H_\infty|}{6}
\wedge 1\right)\right).
$$
It is clear that $k(x,t)\downarrow 0$ as  $t\uparrow\infty$.
Thus as $t\to\infty$, $e^{\lambda_1t}T_t(k(\cdot,t))(x)\to 0,$
which implies
\begin{equation}\label{8.27}
  \lim_{t\to\infty}\langle R(-e^{\lambda_1 t/2}\theta_2,\cdot),X_t\rangle= 0\quad \mbox{in probability}.
\end{equation}
Since $V\in L^2(E,m)\cap L^4(E,m)$, by Remark \ref{rem:large}, we have
\begin{equation}\label{8.13}
\lim_{t\to\infty} e^{\lambda_1 t}\langle V , X_t\rangle=
\langle V,~\phi_1\rangle_m W_\infty \quad \mbox{in probability}.
\end{equation}
Therefore, combining \eqref{8.27} and \eqref{8.13}, we get
\begin{equation}\label{8.23}
  \lim_{t\to\infty}\exp\left\{\langle J(t,\cdot),X_t\rangle\right\}
 =\exp\{-\frac{1}{2}\theta_2^2\langle V,~\phi_1\rangle_m W_\infty\}\quad \mbox{ in probability}.
\end{equation}
Since the real part of $J(t,x)$ is less than $0$,
\begin{equation}\label{8.8}
  \left|\exp\left\{\langle J(t,\cdot),X_t\rangle\right\}\right|\le 1.
\end{equation}
Recall that
$\lim_{t\to\infty}e^{\lambda_1 t}\langle\phi_1,X_t\rangle=W_\infty,$ $\P_\mu$-a.s.
Thus by \eqref{8.23}, \eqref{8.8} and the dominated convergence theorem, we get that as $t\to\infty$,
\begin{eqnarray}\label{8.14}
&&\left|\P_{\mu}\left(\exp\left\{\left(i\theta_1
 -\frac{1}{2}\theta_2^2\langle V,~\phi_1\rangle_m\right) e^{\lambda_1t}\langle \phi_1,X_t\rangle
  +i\theta_3t^{-1/2}e^{\lambda_1 t/2}\langle h,X_t\rangle
   +i\theta_4e^{\lambda_1 t/2}\langle f,X_t\rangle\right\}\right)\right.\nonumber\\
   &&\quad -\kappa_1(\theta_1,\theta_2,\theta_3,\theta_4,t)|\nonumber\\
  &\le&\P_{\mu}\left|\exp\left\{\langle J(t,\cdot),X_t\rangle\right\}
  -\exp\left\{-\frac{1}{2}\theta_2^2\langle V,~\phi_1\rangle_m
  e^{\lambda_1t}\langle \phi_1,X_t\rangle \right\}\right|
  \to 0.
\end{eqnarray}
By Lemma \ref{lem:cs},
 \begin{eqnarray}\label{8.26}
  &&\lim_{t\to\infty}\P_{\mu}\left(\exp\left\{\left(i\theta_1
 -\frac{1}{2}\theta_2^2\langle V,~\phi_1\rangle_m\right) e^{\lambda_1t}\langle \phi_1,X_t\rangle
  +i\theta_3t^{-1/2}e^{\lambda_1 t/2}\langle h,X_t\rangle
   +i\theta_4e^{\lambda_1 t/2}\langle f,X_t\rangle\right\}\right)\nonumber\\
  &=&\P_{\mu}\left(\exp\{i\theta_1W_\infty\}
  \exp\left\{-\frac{1}{2}(\theta_2^2\langle V,~\phi_1\rangle_m
  +\theta_3^2\rho^2_f+\theta_4^2\sigma^2_f) W_\infty\right\}\right).
\end{eqnarray}
By \eqref{L2H}, we get
\begin{equation*}
\langle V,~\phi_1\rangle_m
=\int_0^\infty e^{-\lambda_1s}\left\langle A(I_sg)^2,\phi_1\right\rangle_m\,ds.
\end{equation*}
The proof is now complete.\hfill$\Box$

\smallskip

{\bf Acknowledgements.} We thank the referee for helpful comments on the
first version of this paper.

\smallskip

\begin{singlespace}

\end{singlespace}

\vskip 0.2truein
\vskip 0.2truein

\noindent{\bf Yan-Xia Ren:} LMAM School of Mathematical Sciences \& Center for
Statistical Science, Peking
University,  Beijing, 100871, P.R. China. Email: {\texttt
yxren@math.pku.edu.cn}

\smallskip
\noindent {\bf Renming Song:} Department of Mathematics,
University of Illinois,
Urbana, IL 61801, U.S.A.
Email: {\texttt rsong@math.uiuc.edu}

\smallskip

\noindent{\bf Rui Zhang:} LMAM School of Mathematical Sciences, Peking
University,  Beijing, 100871, P.R. China. Email: {\texttt
ruizhang8197@gmail.com}

\end{doublespace}


\begin{thebibliography}{99}

\bibitem{RP} Adamczak, R. and Mi{\l}o\'{s}, P.:
CLT for Ornstein-Uhlenbeck branching particle system. Preprint, 2011.
arXiv:1111.4559.

\bibitem{AH83} Asmussen, S. and Hering, H.: \emph{Branching Processes}.
Birkh\"auser, Boston, 1983.

\bibitem{AK} Asmussen, S. and Keiding, N.:
Martingale central limit theorems and asymptotic estimation theory
for multitype branching processes. \emph{Adv. Appl. Probab.} \textbf{10}
(1978), 109--129.

\bibitem{Ath69a} Athreya, K. B.:
Limit theorems for multitype continuous time Markov branching processes
I: The case of an eigenvector linear functional.
\emph{Z. Wahrs. Verw. Gebiete} \textbf{12} (1969), 320--332.


\bibitem{Ath69}Athreya,  K. B.:
Limit theorems for multitype continuous time Markov branching processes
II: The case of an arbitrary linear functional.
\emph{Z. Wahrs. Verw. Gebiete} \textbf{13} (1969), 204--214.


\bibitem{Ath71} Athreya, K. B.:
Some refinements in the theory of supercritical multitype Markov branching processes.
\emph{Z. Wahrs. Verw. Gebiete}  \textbf{20} (1971), 47--57.

\bibitem{DS} Davies, E. B. and Simon, B.: Ultracontractivity and the
kernel for Schr\"{o}dinger operators and Dirichlet Laplacians.
\emph{J. Funct. Anal.} \textbf{59} (1984), 335--395.


\bibitem{Dawson} Dawson, D. A. : \emph{Measure-Valued Markov Processes}. Springer-Verlag, 1993.

\bibitem{Dudley} Dudley, R. M.:
\emph{Real Analysis and Probability}. Cambridge University Press, 2002.


\bibitem{E.B.} Dynkin, E. B.: Superprocesses and partial differential equations. \emph{Ann. Probab.}
\textbf{21} (1993), 1185--1262.


\bibitem{E.B2.} Dynkin, E. B. and Kuznetsov, S. E.:
$\mathbb N$-measure for branching exit Markov system and  their applications to differential equations.
\emph{Probab. Theory Rel. Fields} \textbf{130} (2004), 135--150.

\bibitem{elk-roe} El Karoui, N. and Roelly, S.:
Propri\'et\'es de martingales, explosion et repr\'esentation de
L\'evy-Khintchine d'une classe de processus de branchment \`a valeurs mesures.
\emph{Stoch. Proc. Appl.} \textbf{38} (1991), 239--266.


\bibitem{KS} Kesten, H. and Stigum, B. P.:
A limit theorem for multidimensional Galton-Watson processes.
\emph{Ann. Math. Statist.} \textbf{37} (1966), 1211--1223.

\bibitem{KS66} Kesten, H. and Stigum, B. P.:
Additional limit theorems for indecomposable multidimensional Galton-Watson processes.
\emph{Ann. Math. Statist.} \textbf{37} (1966), 1463--1481.


\bibitem{KPR} Kyprianou,  A. E., Perez, J-L. and Ren, Y.-X.:
The backbone decomposition for spatially dependent supercritical superprocesses.
To appear in \emph{S\'{e}minaire de Probabilit\'{e}. }

\bibitem{Li03} Li, Z.: Skew convolution semigroups and related immigration
processes. \emph{Theory Probab. Appl.} \textbf{46} (2003), 274--296.

\bibitem{Li11} Li, Z.: \emph{Measure-valued Branching Markov Processes}. Springer,
Heidelberg, 2011.


\bibitem{Mi} Mi{\l}o\'{s}, P.: Spatial CLT for the supercritical  Ornstein-Uhlenbeck
superprocess. Preprint, 2012. arXiv:1203:6661

\bibitem{RSZ} Ren, Y.-X., Song, R. and Zhang, R.:
Central limit theorems for super Ornstein-Uhlenbeck processes,
\emph{Acta Appl. Math.} \textbf{130} (2014), 9--49.

\bibitem{RSZ2} Ren, Y.-X., Song, R. and Zhang, R.:
Central limit theorems for supercritical branching Markov processes.
\emph{J. Funct. Anal.} \textbf{266} (2014), 1716--1756.

\end{thebibliography}
\end{document}